\documentclass{article}

\usepackage{amsthm,amssymb,amsmath,enumerate,graphicx,psfrag}
\usepackage{algorithm}
\usepackage{algorithmic}
\usepackage{fullpage}
\usepackage{url}
\usepackage{authblk}
\usepackage{caption}

\usepackage{hyperref}

\newtheorem{definition}{Definition}
\newtheorem{claim}{Claim}

\newtheorem{theorem}[definition]{Theorem}

\newtheorem{lemma}[definition]{Lemma}

\usepackage{color}

\newcommand{\comment}[1]{}

\title{Obstructions for three-coloring graphs without induced paths on six vertices}

\author[1]{Maria Chudnovsky\footnote{Partially supported
by NSF grant DMS-1265803.}}
\author[2]{Jan Goedgebeur\footnote{Supported by a Postdoctoral Fellowship of the Research Foundation Flanders (FWO).}}
\author[3]{Oliver Schaudt}
\author[4]{Mingxian Zhong}

\affil[1]{\textit{\small Princeton University, Princeton, NJ 08544,
USA. E-mail: mchudnov@math.princeton.edu}}
\affil[2]{\textit{\small Ghent University, Ghent, Belgium. E-mail: jan.goedgebeur@ugent.be}}
\affil[3]{\textit{\small Universit\"at
zu K\"oln, K\"oln, Germany. E-mail: schaudto@uni-koeln.de}}
\affil[4]{\textit{\small Columbia University, New York, NY 10027,
USA. E-mail: mz2325@columbia.edu}}

\date{}

\begin{document}

\maketitle

\begin{abstract}
We prove that there are 24 4-critical $P_6$-free graphs, and give the complete list.
We remark that, if $H$ is connected and not a subgraph of $P_6$, there are infinitely many 4-critical $H$-free graphs.
Our result answers questions of Golovach et al.~and Seymour.
\end{abstract}

\section{Introduction}

A \emph{$k$-coloring} of a graph $G=(V,E)$ is a mapping $c : V \to \{1,\ldots,k\}$ such that $c(u) \neq c(v)$ for all edges $uv \in E$.
If a $k$-coloring exists, we say that $G$ is \emph{$k$-colorable}.
We say that $G$ is \emph{$k$-chromatic} if it is $k$-colorable but not $(k-1)$-colorable.
A graph is called \emph{$k$-critical} if it is $k$-chromatic, but every proper subgraph is $(k-1)$-colorable.
For example, the class of $3$-critical graphs is the family of all chordless odd cycles.
The characterization of critical graphs is a notorious problem in the theory of graph coloring, and also the topic of this paper.

Since it is NP-hard to decide whether a given graph admits a $k$-coloring, assuming $k \ge 3$, there is little hope of giving a characterization of the $(k+1)$-critical graphs that is useful for algorithmic purposes.
The picture changes if one restricts the structure of the graphs under consideration.

Let a graph $H$ and a number $k$ be given.
An \emph{$H$-free} graph is a graph that does not contain $H$ as an induced subgraph.
We say that a graph $G$ is \emph{$k$-critical $H$-free} if $G$ is $H$-free, $k$-chromatic, and every $H$-free proper subgraph of $G$ is $(k-1)$-colorable.
Note that a \emph{$k$-critical $H$-free} graph is not necessarily a $k$-critical graph.
In this paper we stick to the case of $4$-critical graphs; these graphs we may informally call \emph{obstructions}.

\begin{figure}[t!]
	\centering
		\includegraphics[width=0.7\textwidth]{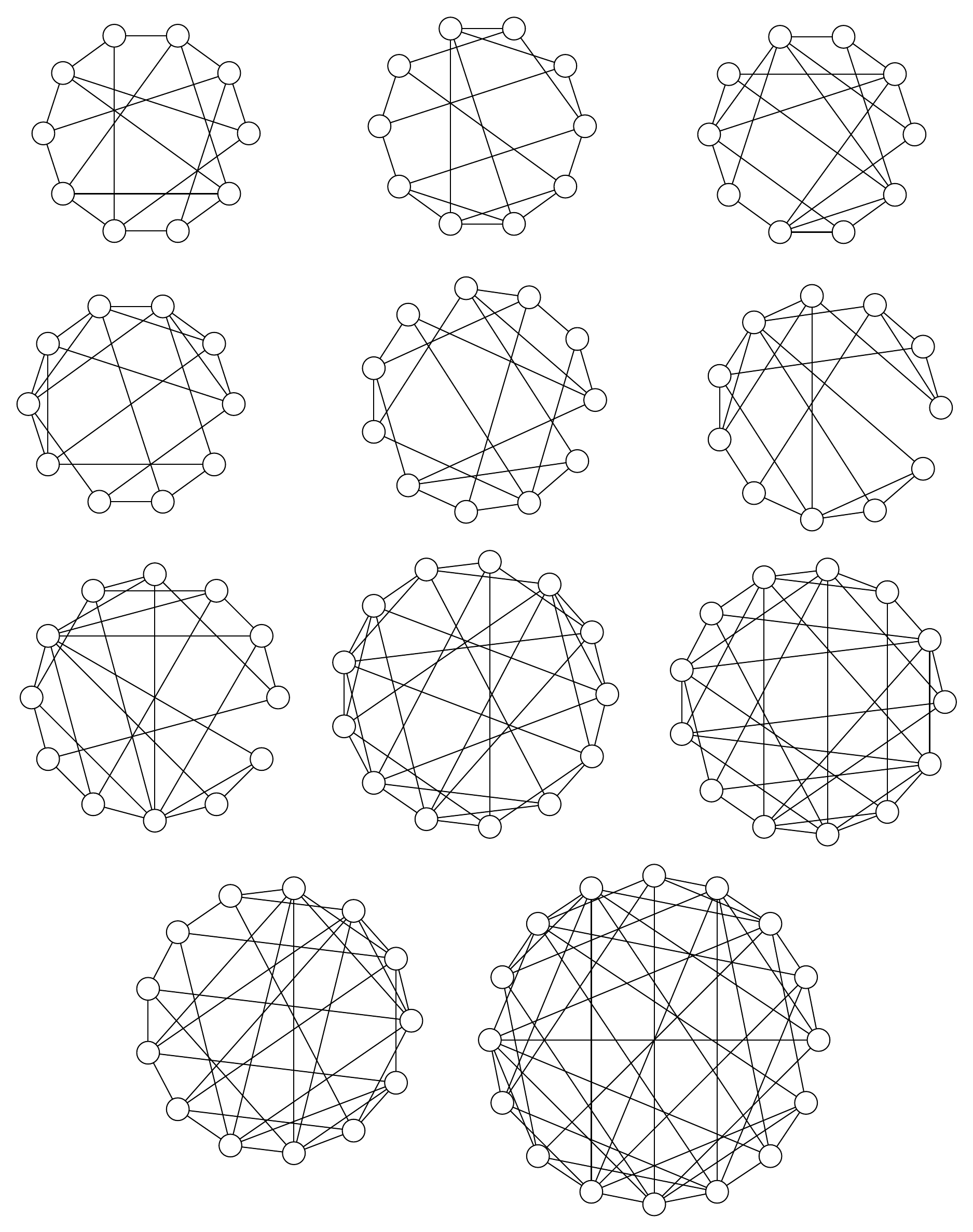}
	\caption{The graphs $F_{14}$ to $F_{24}$, in reading direction}
	\label{fig:graphs14-24}
\end{figure}

Bruce et al.~\cite{BHS09} proved that there are exactly six 4-critical $P_5$-free graphs, where $P_t$ denotes the path on $t$ vertices.
Randerath et al.~\cite{RST02} have shown that the only 4-critical $P_6$-free graph without a triangle is the Gr\"otzsch graph (i.e., the graph $F_{18}$ in Fig.~\ref{fig:graphs14-24}).
More recently, Hell and Huang~\cite{hell_14} proved that there are four 4-critical $P_6$-free graphs without induced four-cycles.

In view of these results, Golovach et al.~\cite{GJPS15} posed the question of whether the list of 4-critical $P_6$-free graphs is finite (cf.~\emph{Open Problem 4} in~\cite{GJPS15}).
In fact, they ask whether there is a certifying algorithm for the 3-colorability problem in the class of $P_6$-free graphs, which is an immediate consequence of the finiteness of the list.
Our main result answers this question affirmatively.

\begin{theorem}\label{thm:N3P6}
There are exactly 24 $4$-critical $P_6$-free graphs.
\end{theorem}

\begin{figure}[t!]
	\centering
		\includegraphics[width=0.7\textwidth]{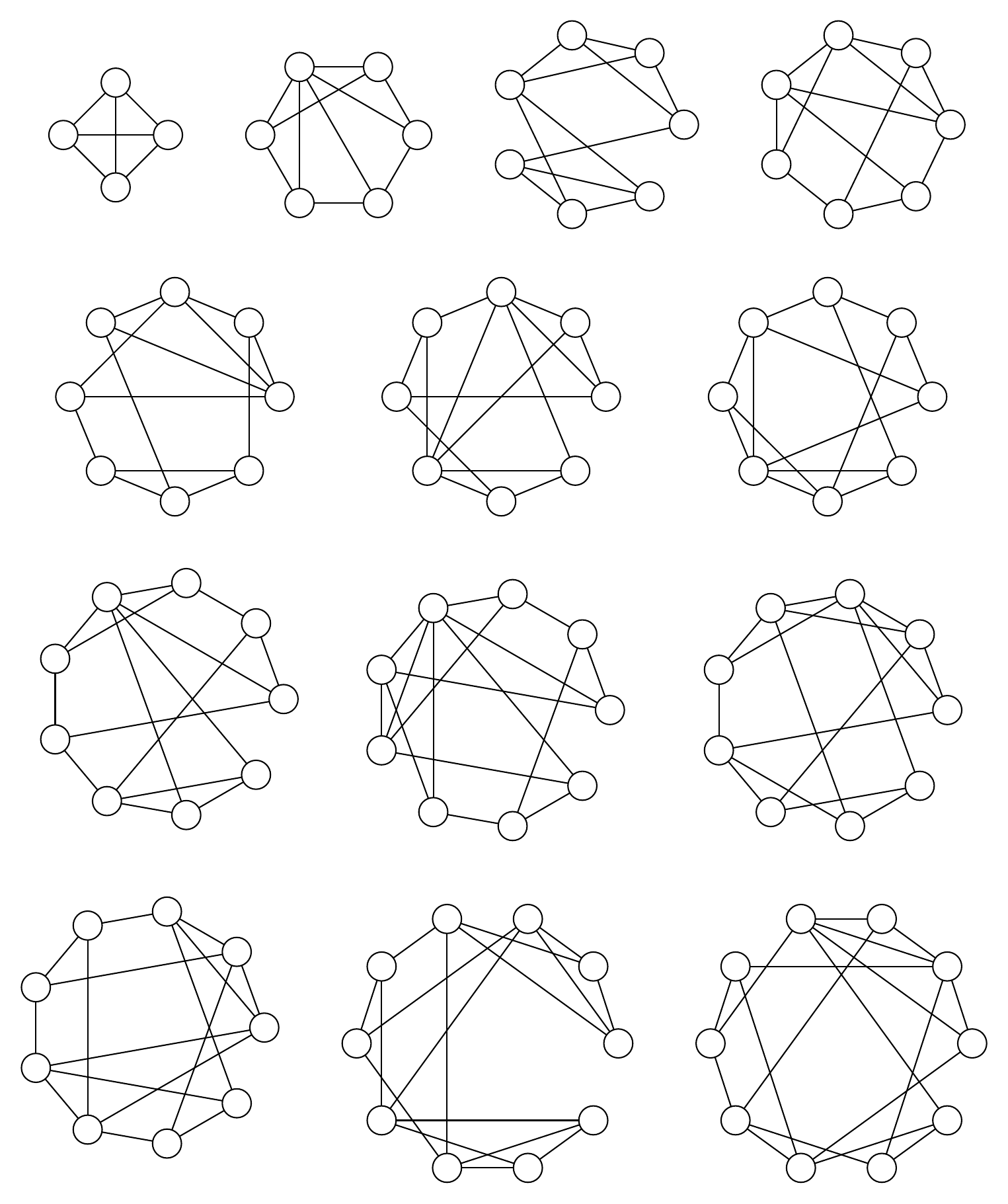}
	\caption{The graphs $F_1$ to $F_{13}$, in reading direction}
	\label{fig:graphs1-13}
\end{figure}

These 24 graphs, which we denote here by $F_1$-$F_{24}$, are shown in Fig.~\ref{fig:graphs1-13} and~\ref{fig:graphs14-24}.
The list contains several familiar graphs, e.g., $F_1$ is $K_4$, $F_2$ is the 5-wheel, $F_3$ is the Moser-spindle, and $F_{18}$ is the Gr\"otzsch graph.
The adjacency lists of these graphs can be found in the Appendix.
Let $\mathcal L = \{F_1,\ldots,F_{24}\}$.

We also determined that there are exactly 80 4-vertex-critical $P_6$-free graphs (details on how we obtained these graphs can be found in the Appendix). 
Table~\ref{table:counts_animals} gives an overview of the counts of all 4-critical and 4-vertex-critical $P_6$-free graphs. All of these graphs can also be obtained from the \textit{House of Graphs}~\cite{hog} by searching for the keywords ``4-critical P6-free'' or ``4-vertex-critical P6-free'' where several of their invariants can be found.  

In Section~\ref{sec:P7} we show that there are infinitely many $4$-critical $P_7$-free graphs using a construction due to Pokrovskiy~\cite{Pok14}.
Note that there are infinitely many 4-critical claw-free graphs.
For example, this follows from the existence of 4-regular bipartite graphs of arbitrary large girth (cf.~\cite{LU95} for an explicit construction of these), whose line graphs are then 4-chromatic.
Also, there are 4-chromatic graphs of arbitrary large girth, which follows from a classical result of Erd\H{o}s ~\cite{Erd59}.
This together with Theorem~\ref{thm:N3P6} yields the following dichotomy theorem, which answers a question of Seymour~\cite{SeyPriv}.

\begin{theorem}
Let $H$ be a connected graph. 
There are finitely many $4$-critical $H$-free graphs if and only if $H$ is a subgraph of $P_6$. 
\end{theorem}

We will next give a sketch of the proof of our main result, thereby explaining the structure of this paper.
An extended abstract of this paper appeared at SODA 2016~\cite{DBLP:conf/soda/ChudnovskyGSZ16}.

\begin{table}[ht!]
\centering
\begin{tabular}{| c | r | r |}
\hline 
Vertices & Critical graphs & Vertex-critical graphs\\
\hline 
4  &  1  &  1\\
6  &  1  &  1\\
7  &  2  &  7\\
8  &  3  &  6\\
9  &  4  &  16\\
10  &  6  &  34\\
11  &  2  &  3\\
12  &  1  &  1\\
13  &  3  &  9\\
16  &  1  &  2\\
\hline
total  &  24  &  80\\
\hline
\end{tabular}
\caption{Counts of all 4-critical and 4-vertex-critical $P_6$-free graphs}

\label{table:counts_animals}
\end{table}

\subsection{Sketch of the proof}

Given a 4-critical $P_6$-free graph $G$, our aim is to show that it is contained in $\mathcal L$.
Our proof is based on the contraction (and uncontraction) of a particular kind of subgraph called \emph{tripod}.
Tripods have been used before in the design of 3-coloring algorithms for $P_7$-free graphs~\cite{BCMSSZ15}.
In Section~\ref{sec:tripods} tripods are defined, and it is shown that contracting a maximal tripod to a single triangle is a \textit{safe} operation for our purpose.

Our proof considers two cases.

\begin{enumerate}
\item If all maximal tripods of $G$ are just single triangles, $G$ is diamond-free. Here, a \emph{diamond} is the graph obtained by removing an edge from $K_4$. 
We determine all 4-critical $(P_6,\mbox{diamond})$-free graphs in Section~\ref{sec:diamond-free} and verify that they are in $\mathcal L$.
Our proof is computer-aided and builds on a substantial strengthening of a method by Ho\`ang et al.~\cite{HMRSV15}.

\item If $G$ contains a maximal tripod that is not a triangle, the structural analysis is more involved. We further distinguish two cases.
\begin{enumerate}[(a)]
	\item If $G$ has a vertex $x$ with neighbors in all three sets of the tripod, an exhaustive computer search shows that $G$ has at most 18 vertices.
	This is done in Section~\ref{sec:tripod-generator}.
The algorithm used is completely different from the one used in the $(P_6,\mbox{diamond})$-free case.
It mimics the way that a tripod can be traversed, thereby applying a set of strong pruning rules that exploit the minimality of the obstruction.

In Section~\ref{sec:smallgraphs} we show that a $4$-critical $P_6$-free graph on at most 28 vertices is a member of $\mathcal L$ using a computer search similar to the one in Section~\ref{sec:diamond-free}.
This time, we allow the graphs to contain diamond, but discard it if it has more than 28 vertices.
From this result we see that $G$ must be a member of $\mathcal L$ since $|V(G)|\le 18$.

\item
If $G$ does not contain such a vertex, a rigorous structural analysis (presented in Section~\ref{sec:uncontraction}) combined with an inductive argument using the previous results proves that $|V(G)|\le 28$.
Like in case (a), it follows that $G \in \mathcal L$.

The analysis in Section~\ref{sec:uncontraction} does not rely on a computer search.
\end{enumerate}
\end{enumerate}

We wrap up the whole proof in Section~\ref{sec:proof-of-main-result}.

As mentioned earlier, in Section~\ref{sec:P7} we show that there are infinitely many $4$-critical $P_7$-free graphs, which results in our dichotomy theorem.

\section{Tripods}\label{sec:tripods}

Let $G$ be a graph, and let $A \subseteq V(G)$ and $b \in V(G) \setminus A$. 
We say that $b$ is {\em complete} to $A$ if $b$ is adjacent to every vertex of $A$, and $b$ is {\em anticomplete} to $A$ if $b$ is non-adjacent to every vertex of $A$. If $b$ has both a neighbor and a non-neighbor in $A$, then
$b$ is {\em mixed on} $A$.
For $B \subseteq V(G) \setminus A$, $B$ is {\em complete} to $A$
if every vertex of $B$ is complete to $A$, and $B$ is {\em anticomplete} to 
$A$ if every vertex of $B$ is anticomplete to $A$.

A \emph{tripod} in a graph $G$ is a triple $T=(A_1,A_2,A_3)$ of disjoint stable sets with the following properties:
\begin{enumerate}[(a)]
	\item $A_1 \cup A_2 \cup A_3 = \{v_1,\ldots,v_k\}$;
	\item $v_i \in A_i$ for $i=1,2,3$;
	\item $v_1v_2v_3$ is a triangle, the \emph{root} of $T$; and
	\item for all $i \in \{1,2,3\}$, $\{\ell,k\}=\{1,2,3\}\setminus\{i\}$, and $j \in \{4,\ldots,k\}$ with $v_j \in A_i$, the vertex $v_j$ has neighbors in both $\{v_1,\ldots,v_{j-1}\} \cap A_{\ell}$ and $\{v_1,\ldots,v_{j-1}\} \cap A_k$.
\end{enumerate}
Assuming that $G$ admits a 3-coloring, it follows right from the definition above that each $A_i$ is contained in a single color class.
Moreover, since $v_1v_2v_3$ is a triangle, $A_1,A_2,A_3$ are pairwise contained in distinct color classes.

\begin{figure}
	\centering
		\includegraphics[width=0.60\textwidth]{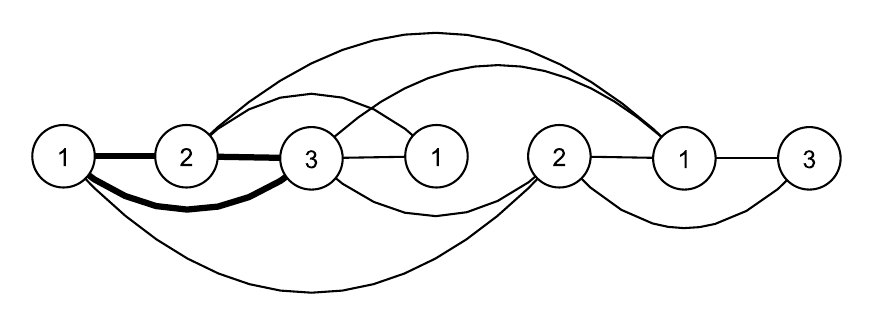}
	\caption{A tripod $T$ with vertex set $V(T)=\{v_1,\ldots,v_7\}$ from left to right. The bold triangle is the root of $T$. The numbers in the vertices denote the unique coloring (up to permutation of colors) of $T$.}
	\label{fig:tripod-illustration}
\end{figure}

An illustration of a tripod is given in Fig.~\ref{fig:tripod-illustration}.

To better reference the ordering of the tripod, we put $t(v_1)=t(v_2)=t(v_3)=0$, and $t(v_i)=i-3$ for all $4 \le i \le k$.
For each $u \in A_i$, let $n_j(u)$ be the neighbor $v$ of $u$ in $A_j$ with $t(v)$ minimum, where $i,j \in \{1,2,3\}, i \neq j$.   
We write $T(t) = G|\{v \in V(T) : t(v) \le t\}$, i.e., the subgraph induced by $G$ on the vertex set $\{v \in V(T) : t(v) \le t\}$.
Moreover, we write $T_i$ for the graph $G|(A_j \cup A_k)$ where $\{i,j,k\} = \{1,2,3\}$, and finally $T_i(t)$ for the graph $G|\{v \in A_j \cup A_k : t(v) \le t\}$.

We call a tripod $(A_1,A_2,A_3)$ \emph{maximal} in a given graph if no further vertex can be added to any set $A_i$ without violating the tripod property.

\subsection{Contracting a tripod}

By \emph{contracting} a tripod $(A_1,A_2,A_3)$ we mean the operation of identifying each $A_i$ to a single vertex $a_i$, for all $i = 1,2,3$.
We then make $a_i$ adjacent to the union of neighbors of the vertices in $A_i$, for all $i = 1,2,3$.

The neighborhood of a vertex $v$ in a graph $G$ is denoted $N_G(v)$.
If $G$ is clear from the context we might also omit $G$ in the subscript.
If $U$ is a subset of the vertex set of $G$, we denote the induced subgraph on the set $U$ by $G|U$.

\begin{lemma}\label{lem:contraction-is-safe}
Let $G$ be a graph with a maximal tripod $T$ such that no vertex of $G$ has neighbors in all three classes of $T$.
Let $G'$ be the graph obtained from $G$ by contracting $T$.
Then the following holds.
\begin{enumerate}[(a)]
	\item\label{lem:contraction-is-safe.coloring} The graph $G$ is 3-colorable if and only if $G'$ is 3-colorable, and
	\item\label{lem:contraction-is-safe.P6-free} if $G$ is $P_6$-free, $G'$ is $P_6$-free.
\end{enumerate}
\end{lemma}
\begin{proof}
Assertion~\eqref{lem:contraction-is-safe.coloring} follows readily from the definition of a tripod, so we just prove~\eqref{lem:contraction-is-safe.P6-free}.
For this, suppose that $G$ is $P_6$-free but $G'$ contains an induced $P_6$, say $P = v_1$-$\ldots$-$v_6$.
Let $T=(A_1,A_2,A_3)$, and let $a_i$ be the vertex of $G'$ the set $A_i$ is contracted to, for $i=1,2,3$.

Since $P$ is an induced path, it cannot contain all three of $a_1,a_2,a_3$.
Moreover, if $P$ contains neither of $a_1,a_2,a_3$, then $G$ contains a $P_6$, a contradiction.

Suppose that $P$ contains, say, $a_1$ and $a_2$.
We may assume that $a_1 = v_i$ and $a_2 =v_{i+1}$ for some $1 \le i \le 3$.
If $i=1$, pick $b \in A_1$ and $c \in A_2 \cap N_G(v_3)$ with minimum distance in $T_3$.
Otherwise, if $i \ge 2$, pick $b \in A_1 \cap N_G(v_{i-1})$ and $c \in A_2 \cap N_G(v_{i+2})$ again with minimum distance in $T_3$.
In both cases, let $Q$ be the shortest path between $b$ and $c$ in $T_3$.
Such a path exists since $T_i$ is connected for $i=1,2,3$ by the definition of a tripod.
Due to the choice of $b$ and $c$, the induced path $v_1$-$\ldots$-$b$-$Q$-$c$-$\ldots$-$v_6$ is induced in $G$, which means $G$ contains a $P_6$, a contradiction.

So, we may assume that $P$ contains only one of $a_1,a_2,a_3$, say $v_i=a_1$ for some $1 \le i \le 3$.
We obtain an immediate contradiction if $i=1$, so suppose that $i \ge 2$.
Since $v_{i+2}$ is not contained in $T$, we may assume that~$v_{i+2}$ is anticomplete to $A_2$ in $G$.
Pick $b \in A_1 \cap N_G(v_{i-1})$ and $c \in A_1 \cap N_G(v_{i+1})$ such that the distance in $T_3$ between $b$ and $c$ is minimum.
Let $Q$ be a shortest path in $T_3$ between $b$ and $c$.
Since $v_iv_{i+2} \notin E(G')$, $v_{i+2}$ is anticomplete to $A_1$ and thus to $V(Q)$ in $G$.
If $b=c$, then $v_1$-$\ldots$-$v_{i-1}$-$b$-$v_{i+1}$-$\ldots$-$v_6$ is a $P_6$ in $G$, a contradiction.
Otherwise, the induced path $v_{i-1}$-$b$-$Q$-$c$-$v_{i+1}$-$v_{i+2}$ is induced in $G$ and contains at least six vertices, which is also contradictory.
\end{proof}

\section{Diamond-free obstructions}\label{sec:diamond-free}

Recall that a diamond is the graph obtained by removing an edge from $K_4$.
After successively contracting all maximal tripods in a graph, we are left with a diamond-free graph.
In this section we prove the following statement.

\begin{lemma}\label{lem:diamondfree}
There are exactly six 4-critical $(P_6,\mbox{diamond})$-free graphs.
\end{lemma}

These graphs are $F_1$, $F_{11}$, $F_{14}$, $F_{16}$, $F_{18}$, and $F_{24}$ in Fig.~\ref{fig:graphs1-13} and~\ref{fig:graphs14-24}.

The proof of Lemma~\ref{lem:diamondfree} is computer-aided and builds upon a method recently proposed by Ho\`ang et al.~\cite{HMRSV15}.
With this method they have shown that there is a finite number of 5-critical $(P_5,C_5)$-free graphs. 
The idea is to automatize the large number of necessary case distinctions, resulting in an exhaustive enumeration algorithm.
Since we have to deal with a graph class which is substantially less structured, we need to significantly extend their method.

\subsection{Preparation}

In order to prove Lemma~\ref{lem:diamondfree}, we make use of the following tools.

Let $G$ be a $k$-colorable graph.
We define the \emph{$k$-hull} of $G$, denoted $G_k$, to be the graph with vertex set $V(G)$ where two vertices $u,v$ are adjacent if and only if there is no $k$-coloring of $G$ where $u$ and $v$ receive the same color.
Note that $G_k$ is a simple supergraph of $G$, since adjacent vertices can never receive the same color in any coloring.
Moreover, $G_k$ is $k$-colorable.

It is easy to see that a $k$-critical graph cannot contain two distinct vertices, $u$ and $v$ say, such that $N(u) \subseteq N(v)$.
The following observation is a proper generalization of this fact.
If $U$ is a vertex subset of a graph $G$, we denote by $G-U$ the graph obtained from $G$ by deleting the vertices in $U$.

\begin{lemma}\label{lem:color-trick}
Let $G=(V,E)$ be a $k$-vertex-critical graph and let $U,W$ be two non-empty disjoint vertex subsets of $G$.
Let $H:=(G-U)_{k-1}$.
If there exists a homomorphism $\phi : G|U \mapsto H|W$, then $N_G(u)\setminus U \not\subseteq N_H(\phi(u))$ for some $u \in U$.
\end{lemma}

Note that, in the statement of Lemma~\ref{lem:color-trick}, $H$ is well-defined since $G$ is $k$-vertex-critical.

\begin{proof}[Proof of Lemma~\ref{lem:color-trick}]
Suppose that $N_G(u)\setminus U \subseteq N_H(\phi(u))$ for all $u \in U$.
Fix some $(k-1)$-coloring $c$ of $H$.
In particular, for each $u \in U$, the color of $\phi(u)$ is different from that of any member of $N_H(\phi(u))$.

We now extend $c$ to a $(k-1)$-coloring of $G$ by giving any $u \in U$ the color $c(\phi(u))$.
It suffices to show that this is a proper coloring.
Clearly there are no conflicts between any two vertices of $U$, since $\phi$ is a homomorphism.
Let $u \in U$ and $v \in N_G(u)\setminus U$ be arbitrary.
Since $N_G(u)\setminus U \subseteq N_H(\phi(u))$, $c(v) \neq c(\phi(u))$, and so $u$ and $v$ receive distinct colors.
But this contradicts with the assumption that $G$ is a $k$-vertex-critical graph.
\end{proof}

We make use of Lemma~\ref{lem:color-trick} in the following way.
Assume that $G$ is a $(k-1)$-colorable graph that is an induced subgraph of some $k$-vertex-critical graph $G'$.
Pick two non-empty disjoint vertex subsets $U,W \subseteq V$ of $G$, and let $H:=(G-U)_{k-1}$.
Assume there exists a homomorphism $\phi : G|U \mapsto H|W$ such that $N_G(u)\setminus U \subseteq N_H(\phi(u))$ for all $u \in U$.
Then there must be some vertex $x \in V(G')\setminus V(G)$ which is adjacent to some $u \in U$ but non-adjacent to $\phi(u)$ in $G'$.
Moreover, $x$ is non-adjacent to $\phi(u)$ in the graph $(G'-U)_{k-1}$.

We also make use of the following well-known fact.

\begin{lemma}\label{lem:min_degree}
A $k$-vertex-critical graph has minimum degree at least $k-1$.
\end{lemma}

Another fact we need is the following.

\begin{lemma}\label{lem:C5-must-be-there}
Any $(P_6,\mbox{diamond})$-free 4-critical graph other than $K_4$ contains an induced $C_5$.
\end{lemma}
\begin{proof}
By the Strong Perfect Graph Theorem~\cite{chudnovsky_06}, every $4$-critical graph different from $K_4$ must contain an odd hole or an anti-hole as an induced subgraph.
A straightforward argumentation shows that only the 5-hole, $C_5$, can possibly appear.
\end{proof}

\subsection{The enumeration algorithm}

Generally speaking, our algorithm (Algorithm~\ref{algo:init-algo} together with Algorithm~\ref{algo:construct}) constructs a graph $G'$ with $n+1$ vertices from a graph $G$ with $n$ vertices by adding a new vertex and connecting it to vertices of $G$ in all possible ways. So, all graphs constructed from $G$ contain $G$ as an induced subgraph.
Since 3-colorability and $(P_6,\mbox{diamond})$-freeness are both hereditary properties, we do not need to expand $G$ if it is not 3-colorable, contains a $P_6$ or a diamond.

We use Algorithm~\ref{algo:init-algo} below to enumerate all $(P_6,\mbox{diamond})$-free $4$-critical graphs.
In order to keep things short, we use the following conventions for a graph $G$.
We call a pair $(u,v)$ of distinct vertices for which $N_G(u) \subseteq N_{(G-u)_3}(v)$ \emph{similar vertices}.
Similarly, we call a 4-tuple $(u,v,u',v')$ of distinct vertices with $uv, u'v' \in E(G)$ such that $N_G(u)\setminus \{v\} \subseteq N_{(G-\{u,v\})_3}(u')$ and $N_G(v)\setminus \{u\} \subseteq N_{(G-\{u,v\})_3}(v')$ \emph{similar edges}.
Finally, we define \emph{similar triangles} in an analogous fashion.

\begin{algorithm}[h]
\caption{Generate $(P_6,\mbox{diamond})$-free $4$-critical graphs}
\label{algo:init-algo}
  \begin{algorithmic}[1]
	\STATE Let $\mathcal F$ be an empty list
	\STATE Add $K_4$ to the list $\mathcal F$
	\STATE Construct($C_5$) \ // i.e.\ perform Algorithm~\ref{algo:construct}
	\STATE \textbf{return} $\mathcal F$
  \end{algorithmic}
\end{algorithm}

\begin{algorithm}[ht!]
\caption{Construct(Graph $G$)}
\label{algo:construct}
  \begin{algorithmic}[1]
		\IF{$G$ is $(P_6, \mbox{diamond})$-free AND not generated before} \label{line:isocheck}
			\IF{$G$ is not $3$-colorable}
				\IF{$G$ is $4$-critical $P_6$-free} \label{line:k-critical}
					\STATE add $G$ to the list $\mathcal F$
				\ENDIF	
			\STATE \textbf{return}
			\ELSE
				\IF{$G$ contains similar vertices $(u,v)$}
					\FOR{every graph $H$ obtained from $G$ by attaching a new vertex $x$ and incident edges in all possible ways, such that $ux \in E(H)$, but $vx \notin E((H-u)_3)$}
						\STATE\label{line:dominatedvertex} Construct($H$)
					\ENDFOR
				\ELSIF{$G$ contains a vertex $u$ of degree at most $2$}	
					\FOR{every graph $H$ obtained from $G$ by attaching a new vertex $x$ and incident edges in all possible ways, such that $ux \in E(H)$}
						\STATE Construct($H$)
					\ENDFOR
				\ELSIF{$G$ contains similar edges $(u,v,u',v')$}\label{line:dominatededge}	
					\FOR{every graph $H$ obtained from $G$ by attaching a new vertex $x$ and incident edges in all possible ways, such that $ux \in E(H)$ and $u'x \notin E((H-\{u,v\})_3)$, or $vx \in E(H)$ and $v'x \notin E((H-\{u,v\})_3)$}
						\STATE Construct($H$)
					\ENDFOR
				\ELSIF{$G$ contains similar triangles $(u,v,w,u',v',w')$}\label{line:dominatedtriangle}
					\FOR{every graph $H$ obtained from $G$ by attaching a new vertex $x$ and incident edges in all possible ways, such that $ux \in E(H)$ and $u'x \notin E((H-\{u,v,w\})_3)$, $vx \in E(H)$ and $v'x \notin E((H-\{u,v,w\})_3)$, or $wx \in E(H)$ and $w'x \notin E((H-\{u,v,w\})_3)$}
						\STATE Construct($H$)
					\ENDFOR						
				\ELSE 
					\FOR{every graph $H$ obtained from $G$ by attaching a new vertex $x$ and incident edges in all possible ways}
						\STATE\label{line:smallvertex} Construct($H$)
					\ENDFOR		
				\ENDIF					
			\ENDIF	
		\ENDIF	  
  \end{algorithmic}
\end{algorithm}

We now prove that Algorithm~\ref{algo:init-algo} is correct.

\begin{lemma}\label{lem:diamond-algo-correct}
Assume that Algorithm~\ref{algo:init-algo} terminates, and outputs the list of graphs $\mathcal F$.
Then $\mathcal F$ is the list of all $(P_6,\mbox{diamond})$-free 4-critical graphs.
\end{lemma}
\begin{proof}
In view of lines~\ref{line:isocheck} and~\ref{line:k-critical} of Algorithm~\ref{algo:construct}, it is clear that all graphs of $\mathcal F$ are 4-critical $(P_6,\mbox{diamond})$-free.
So, it remains to prove that $\mathcal F$ contains all $(P_6,\mbox{diamond})$-free 4-critical graphs.
To see this, we first prove the following claim.

\begin{claim}\label{clm:invariant}
For every $(P_6,\mbox{diamond})$-free $4$-critical graph $F$ other than $K_4$, Algorithm~\ref{algo:construct} applied to $C_5$ generates an induced subgraph of $F$ with $i$ vertices for every $5 \le i \le |V(F)|$.
\end{claim}
We prove this inductively, as an invariant of our algorithm.
Due to Lemma~\ref{lem:C5-must-be-there}, we know that $F$ contains an induced $C_5$, so the claim holds true for $i=5$.

So assume that the claim is true for some $i \ge 5$ with $i < |V(F)|$. 
Let $G$ be the induced subgraph of $F$ with $|V(G)|=i$.
First assume that $G$ contains similar vertices $(u,v)$.
We put $U=\{u\}$, $W=\{v\}$, $H=(F-u)_3$.
Then, by Lemma~\ref{lem:color-trick}, $N_F(u)\setminus U \not\subseteq N_H(v)$.
Hence, there is some vertex $x \in V(F) \setminus V(G)$ which is adjacent to $u$ in $F$, but not to $v$ in $H$.
Following the statement of line~\ref{line:dominatedvertex}, Construct($F|(V(G) \cup \{x\})$) is called.
We omit the discussion of the lines~\ref{line:dominatededge} and~\ref{line:dominatedtriangle}, as they are analogous.

So assume that $G$ contains a vertex $u$ of degree at most 2.
Then, since the minimum degree of any 4-vertex-critical graph is at least 3, there is some vertex $x \in V(F) \setminus V(G)$ adjacent to $u$.
Following the statement of line~\ref{line:smallvertex}, Construct($F|(V(G) \cup \{x\})$) is called.

Finally, if none of the above criteria apply to $G$, the algorithm attaches a new vertex to $G$ in all possible ways, and calls Construct for all of these new graphs.
Since $|V(F)|>|V(G)|$, among these graphs there is some induced subgraph of $F$, and of course this graph has $i+1$ vertices.
This completes the proof of Claim~\ref{clm:invariant}.
\bigskip

Given that the algorithm terminates and $K_4$ is added to the list $\mathcal F$, Claim~\ref{clm:invariant} implies that $\mathcal F$ must contain all 4-critical $(P_6,\mbox{diamond})$-free graphs.
\end{proof}

We implemented this algorithm in C with some further optimizations. To make sure that no isomorphic graphs are accepted (cf.\ line~\ref{line:isocheck} of Algorithm~\ref{algo:construct}), we use the program \verb|nauty|~\cite{nauty-website, mckay_14} to compute a canonical form of the graphs. We maintain a list of the canonical forms of all non-isomorphic graphs which were generated so far and only accept a graph if it was not generated before (and then add its canonical form to the list).

Our program does indeed terminate (in about 2 seconds), and outputs the six graphs $F_1$, $F_{11}$, $F_{14}$, $F_{16}$, $F_{18}$, and $F_{24}$ from Fig.~\ref{fig:graphs1-13} and~\ref{fig:graphs14-24}. Together with Lemma~\ref{lem:diamond-algo-correct} this proves Lemma~\ref{lem:diamondfree}. 
Let us stress the fact that in order for the algorithm to terminate, all proposed expansion rules are needed.

Table~\ref{table:counts_P6_diamondfree} shows the number of non-isomorphic graphs generated by the program. The source code of the program can be downloaded from~\cite{criticalpfree-site} and in the Appendix we describe how we extensively tested the correctness of our implementation.

The second and third author also extended this algorithm which allowed to determine all $k$-critical graphs for several other cases as well (see~\cite{paper-other-cases}).

\begin{table}[ht!]
\centering
\begin{tabular}{| l || c c c c c c c c c c c c |}
\hline 
$|V(G)|$            & 5 & 6 &  7 &  8 &   9 &  10 &  11 &  12 &  13 &  14 &  15 &  16\\
\# graphs generated & 1 & 4 & 16 & 55 & 130 & 230 & 345 & 392 & 395 & 279 & 211 & 170\\
\hline
$|V(G)|$            &  17 & 18 & 19 & 20 & 21 & 22 & 23 & 24 & 25 & 26 & 27 & 28\\
\# graphs generated & 112 & 95 & 74 & 53 & 40 & 32 & 20 & 15 & 12 &  3 &  1 &  0\\
\hline
\end{tabular}
\caption{Counts of the number of non-isomorphic $(P_6,\mbox{diamond})$-free graphs generated by our implementation of Algorithm~\ref{algo:init-algo}}
\label{table:counts_P6_diamondfree}
\end{table}

\section{Uncontracting a triangle to a tripod}\label{sec:uncontraction}

Let $G$ be a $P_6$-free graph that is not 3-colorable.
Assume that $G$ contains a maximal tripod $T=(A_1,A_2,A_3)$ with $A_1 \cup A_2 \cup A_3 = \{v_1,\ldots,v_k\}$.
The aim of this section is to prove the following statement.

\begin{lemma}\label{thm:10}
Let $G'$ be obtained from $G$ by contracting a maximal tripod
$(A_1,A_2,A_3)$ to a triangle $\{a_1,a_2,a_3\}$. Let $H' \in \mathcal L$ be an induced subgraph of $G'$. If $H'=K_4$, assume that
$|V(H') \cap \{a_1,a_2,a_3\}| <3$. Then there exists an induced subgraph $H$ of $G$ that is not 3-colorable with 
at most $|V(H')|+12$ vertices.
\end{lemma}


In Section~\ref{sec:tripod-generator} a similar statement
is proved for the case when $H'=K_4$ and $a_1,a_2,a_3 \in V(H')$. To construct 
the graph $H$ as in Lemma~\ref{thm:10}  we replace each of $a_1,a_2,a_3$ by a 
subset of $A_1,A_2,A_3$ (call these subsets $C_1,C_2,C_3$ respectively) such 
that every vertex of $N_{H'}(a_i) \setminus \{a_1,a_2,a_3\}$ has a neighbor in 
$C_i$. We then add a few more vertices from $A_1 \cup A_2 \cup A_3$, to ensure that in every 3-coloring of $H$ each of the sets $C_1,C_2,C_3$ is monochromatic, and no colors appear in two of them. 
 
The proof of Lemma~\ref{thm:10} is organized as follows.
The first several claims (Claim~\ref{clm:noleaforhat}---Claim~\ref{thm:XYinT}) 
are technical tools we need later. 
In Claim~\ref{thm:G'} we show that for every $i$ we can 
construct $C_i' \subseteq A_i$ and add at most two more vertices so that every 
vertex of $N_{H'}(a_i) \setminus \{a_1,a_2,a_3\}$ has a neighbor in 
$C_i'$, and $C_i'$ is monochromatic in every 3-coloring of the resulting graph.
Claim~\ref{thm:G'} needs a few technical assumptions; 
in Claim~\ref{thm:7} we show that the assumptions of Claim~\ref{thm:G'} 
hold. In Claim~\ref{thm:8} we analyze how the sets $C_i'$ from  
Claim~\ref{thm:G'} interact for two different values of $i$, again, under certain technical assumptions. In Claim~\ref{thm:9} we show that the assumptions of Claim~\ref{thm:8} hold  except in two special cases. Finally, we deal with the two cases not covered by Claim~\ref{thm:8}, and use  Claim~\ref{thm:8} to produce $H$.

We now describe the proof of Lemma~\ref{thm:10} in detail.
Let $C$ be a hole in $G$. 
A \emph{leaf} for $C$ is a vertex $v \in V(G) \setminus V(C)$ with exactly one neighbor in $V(C)$. 
Similarly, a \emph{hat} for $C$ is a vertex in $V(G) \setminus V(C)$ with exactly two neighbors $u,v \in V(C)$, where $u$ is adjacent to $v$.

The following observation is immediate from the fact that $G$ is $P_6$-free.

\begin{claim}\label{clm:noleaforhat}
No $C_6$ in $G$ has a leaf or a hat.
\end{claim}

\begin{claim}\label{Tijconnected}
The graph $T_i(t)$ is connected, for all $i \in \{1,2,3\}$ and $0 \le t \le k$.
\end{claim}
\begin{proof}
This follows readily from the definition of a tripod.
\end{proof}

\begin{claim}\label{thm:2.5}
Let $a \in A_1$, and let $y,z \in V(G) \setminus (A_1 \cup A_2 \cup A_3)$ such that $a$-$y$-$z$ is an induced path, and $z$ is anticomplete to $A_2 \cup A_3$. 
Then $(A_2 \cup A_3) \setminus N(a)$ is stable, and in particular, for $i=2,3$ there exist $n_i \in N(a) \cap A_i$ such that $n_2$ is adjacent to $n_3$.
\end{claim}
\begin{proof}
By the maximality of the tripod, $y$ is anticomplete to $A_2\cup A_3$.
Suppose there are $p_i \in A_i \setminus N(a)$, $i=2,3$, such that $p_2$ is adjacent to $p_3$. 
Since by \ref{Tijconnected} $T_1$ is connected, we can choose $p_2,p_3$ such that, possibly exchanging $A_2$ and $A_3$, $p_2$ has a neighbor $q_3$ in $A_3 \cap N(a)$. 
But now $z$-$y$-$a$-$q_3$-$p_2$-$p_3$ is a $P_6$, a contradiction. 
Since $T_3$ is connected, the second statement of the theorem follows.\end{proof}

A \emph{2-edge matching} are two disjoint edges $ab,cd$ where $ad,cb$ are non-edges.

\begin{claim}\label{XinT}
Let $X$ be a stable set in $V(G)  \setminus (A_1 \cup A_2 \cup A_3)$, such that for every $x,x' \in X$ there exists $p \in V(G) \setminus (A_1 \cup A_2 \cup A_3)$ such that
$p$ is anticomplete to $A_1$ and adjacent to exactly one of $x,x'$.
Assume that there is a 2-edge matching $ax,a'x'$ between $A_1$ and $X$. 
Then 
\begin{enumerate}[(a)]
	\item\label{XinT(a)} there do not exist $n_2 \in A_2$ and $n_3 \in A_3$ such that $\{a,a'\}$ is complete to $\{n_2,n_3\}$, and
	\item\label{XinT(b)} there exists $a'' \in A_1$, with $t(a'') < \max(t(a),t(a'))$ such that $a''$ is complete to $X \cap (N(a) \cup N(a'))$.
\end{enumerate}
\end{claim}
\begin{proof}
Suppose $ax,a'x'$ is such a matching. 
We may assume that $xp$ is an edge. 
Let $P$ be an induced path from $a$ to $a'$ with interior in $A_2 \cup A_3$.
Such a path exists since $T_1$ is connected, and both $a,a'$ have neighbors in $A_2 \cup A_3$. By the maximality of the tripod, $\{a,a'\}$ is anticomplete
to $A_2 \cup A_3$.
If $P$ has at least three edges, then $x$-$a$-$P$-$a'$-$x'$ is a $P_6$, so we may assume that $a,a'$ have a common neighbor $n_2 \in A_2$. 
If $p$ is non-adjacent to $n_2$, then $p$-$x$-$a$-$n_2$-$a'$-$x'$ is a $P_6$, a contradiction. 
So $p$ is adjacent to $n_2$, and therefore $p$ has no neighbor in $A_3$. 
By symmetry, $a,a'$ have no common neighbor in $A_3$, and so \eqref{XinT(a)} follows.

Since $a,a'$ do not have a common neighbor in $A_3$,  there is an induced path $a$-$b$-$c$-$d$-$a'$  from $a$ to $a'$ in $T_2$. 
It follows from the maximality of the tripod that
$(N(a) \cup N(a')) \cap X$ is anticomplete to $A_2 \cup A_3$.
Since $z$-$a$-$b$-$c$-$d$-$a'$ and $a$-$b$-$c$-$d$-$a'$-$z'$ are not a $P_6$ for any $z \in (N(a) \cap X) \setminus N(a')$, and $z' \in (N(a') \cap X) \setminus N(a)$, we deduce that $c$ is complete to  $((N(a) \cap X) \setminus N(a')) \cup ((N(a')\cap X) \setminus N(a))$. 
We may assume that there exists $x'' \in X \cap N(a) \cap N(a')$ such that $c$ is non-adjacent to $x''$, for otherwise \eqref{XinT(b)} holds. 
Now if $p$ is non-adjacent to $x''$, then $p$-$x$-$c$-$d$-$a'$-$x''$ is a $P_6$, and if $p$ is adjacent to $x''$, then $p$-$x''$-$a$-$b$-$c$-$x'$ is a $P_6$, in both cases a contradiction.
This proves \eqref{XinT(b)}.
\end{proof}

\begin{claim}\label{thm:XYinT}
Let $X,Y$ be two disjoint stable sets in $V(G) \setminus (A_1\cup A_2\cup A_3)$  such that every vertex of $X\cup Y$ has a neighbor in $A_1$.
Moreover, assume that the following assertions hold.
\begin{enumerate}
	\item\label{thm:XYinTasp1} For every $x \in X$ and $y \in Y$, either
		\begin{enumerate}[(i)]
			\item $x$ is adjacent to $y$,
			\item $x$ has a neighbor in $V(G)$ that is anticomplete to $A_1$, or
			\item $y$ has a neighbor in $V(G)$ that is anticomplete to $A_1$.
		\end{enumerate}
	\item\label{thm:XYinTasp2} For every $x,x' \in X$ there exists $p \in V(G) \setminus (A_1\cup A_2\cup A_3)$ such that
		\begin{enumerate}[(i)]
			\item $p$ is anticomplete to $A_1$, and
			\item\label{thm:XYinTasp2(ii)} $p$ is adjacent to exactly one of $x,x'$.
		\end{enumerate}
	\item\label{thm:XYinTasp3} The above assertion holds for $Y$ in an analogous way.
	\item\label{thm:XYinTasp4} Let $u,v \in X\cup Y$ be distinct and non-adjacent. Then $N(u) \setminus A_1$ and $N(v) \setminus A_1$ are incomparable.
\end{enumerate}
Then either
\begin{enumerate}[(a)]
	\item\label{thm:XYinToutcome1} there is a vertex $p \in A_1$ which is complete to $X\cup Y$, or
	\item\label{thm:XYinToutcome2} there exist $c,d \in A_1$, $p \in A_2$ and $q \in A_3$, such that $p$ and $q$ are adjacent, $c$ is complete to $X$, $d$ is complete to $Y$, and $\{c,d\}$ is complete to $\{p,q\}$.
\end{enumerate}
\end{claim}
\begin{proof}
After deleting all vertices of $V(G) \setminus (X \cup Y \cup A_2 \cup A_3)$ with a neighbor in $A_1$ (this does not change the hypotheses or the outcomes), we may assume that no vertex of $V(G) \setminus (A_2\cup A_3\cup X\cup Y)$ has a neighbor in $A_1$.
\begin{equation}\label{eqn:XYinT(1)}
\mbox{\em There exist $a,b \in A_1$ such that $a$ is complete to $X$, and $b$ is complete to $Y$.}
\end{equation}
To see~\eqref{eqn:XYinT(1)}, it is enough to show that $a$ exists, by symmetry. 
So, suppose that such an $a$ does not exist. 
Pick $a \in A_1$ with $N(a) \cap X$ maximal, and note that $a$ is not complete to $X$ by assumption. Let $x' \in X \setminus N(a)$, and let 
$a' \in A_1 \cap N(x')$. By the maximality of $N(a) \cap X$, there exists 
$x \in N(a) \cap X$ such that $a'$ is non-adjacent to $x$.
Now $ax,a'x'$ is a 2-edge matching. But now by Claim~\ref{XinT}\eqref{XinT(b)}, there exists $a'' \in A_1$ complete to $(N(a) \cap X) \cup {x'}$, contrary to the choice of $a$. This proves \eqref{eqn:XYinT(1)}.

We may assume that no vertex of $A_1$ is complete to $X\cup Y$, for otherwise Claim~\ref{thm:XYinT}\eqref{thm:XYinToutcome1} holds.
Moreover, we may assume that there exist $x \in X$, and $y \in Y$ such that $ax,by$ 
is a 2-edge matching.  
We choose $a,b$ with $t(a)+t(b)$ minimum, and subject to that $x$ and $y$ are chosen adjacent if possible.
\begin{equation}\label{eqn:XYinT(2)}
\begin{minipage}[c]{0.8\textwidth}\em
There is no $p \in A_1$, with $t(p) < \max (t(a),t(b))$ such that $p$ is complete to $(X \setminus N(b)) \cup (Y \setminus N(a))$.
\end{minipage}\ignorespacesafterend 
\end{equation}
Suppose such a $p$ exists. 
We may assume that $t(a) > t(b)$, and hence $t(p) < t(a)$.
By the choice of $a$ and $b$, $p$ is not complete to $X$, and so there is a 2-edge matching  between $\{b,p\}$ and $X$. 
Thus by Claim~\ref{XinT}\eqref{XinT(b)}, there exists a vertex $p'$ with $t(p') <  \max(t(b), t(p)) < t(a)$  that is complete to $X$, again contrary to the choice of $a$ and $b$. This proves \eqref{eqn:XYinT(2)}. 
\begin{equation}\label{eqn:XYinT(3)}
\mbox{\em Either $a$ is adjacent to $n_2(b)$, or $b$ is adjacent to $n_2(a)$.}
\end{equation}
Suppose that this is false. 
We may assume that $t(n_2(a))>t(n_2(b))$. 
Let $P$ be an induced path from $n_2(a)$ to $n_2(b)$ in $T_3(t(n_2(a)))$. 
Then $n_2(a)$ is the unique neighbor of $a$ in $P$. 
Since $a$-$n_2(a)$-$P$-$n_2(b)$ is not a $P_6$, we  deduce that $P$ has length two, say $P=n_2(a)$-$p$-$n_2(b)$.  
Moreover, since $x'$-$a$-$n_2(a)$-$p$-$n_2(b)$-$b$ is not a $P_6$ for any $x' \in X \setminus N(b)$, we know that $X \setminus N(b)$ is complete to $p$
(recall that by the maximality of the tripod $x'$ is anticomplete to $A_2$).  
Finally, since $y'$-$b$-$n_2(b)$-$p$-$n_2(a)$-$a$ is not a $P_6$ for any $y' \in Y \setminus N(a)$, $p$ is complete to $Y \setminus N(a)$. 
But since $p \in T_3(t(n_2(a))$, we know that $t(p) < t(a) \leq \max (t(a), t(b))$, contrary to \eqref{eqn:XYinT(2)}. 
This proves \eqref{eqn:XYinT(3)}.

By \eqref{eqn:XYinT(3)} and using the symmetry between $A_2$ and $A_3$, we  deduce that for $i=2,3$ there exists $n_i \in A_i$ such that $\{a,b\}$ is complete to $\{n_2,n_3\}$, and each $n_i$ is the smallest neighbor of one of $a,b$ in $A_i$ w.r.t.~their value of $t$.  
We may assume that $n_2$ is non-adjacent no $n_3$, for otherwise Claim~\ref{thm:XYinT}\eqref{thm:XYinToutcome2} holds.
\begin{equation}\label{eqn:XYinT(4)}
\begin{minipage}[c]{0.8\textwidth}\em
Let $z \in V(G) \setminus (A_1\cup A_2\cup A_3\cup X\cup Y)$ be anticomplete to $A_1$.
Then $z$ is not mixed on any non-edge with one end in $X \setminus N(b)$ and the other in $Y \setminus N(a)$. 
In particular, either $x$ is adjacent to $y$, or some $z \in V(G) \setminus (A_1\cup A_2\cup A_3\cup \{x,y\})$ is complete to $\{x,y\}$ and anticomplete to $A_1$.
\end{minipage}\ignorespacesafterend 
\end{equation}
Suppose $z$ is mixed on a non-edge $x',y'$ with $x' \in X \setminus N(b)$, and $y' \in Y \setminus N(a)$. 
From the maximality of the tripod, we may assume that $z$ is anticomplete to $A_2$.
Now one of the induced paths $z$-$x'$-$a$-$n_2$-$b$-$y'$ and $z$-$y'$-$b$-$n_2$-$a$-$x'$ is a $P_6$, a contradiction.
The second statement follows from assumption \ref{thm:XYinTasp1}. 
This proves \eqref{eqn:XYinT(4)}.

By symmetry, we may assume that $t(n_2) > t(n_3)$, and that $n_2=n_2(a)$.
Thus, there is an induced path $n_2$-$n_3'$-$c$-$n_3$ in $T_1(t(n_2))$. 
Hence $t(c) < t(n_2)$, and so $a$ is non-adjacent to $c$. 
\begin{equation}\label{eqn:XYinT(5)}
\mbox{\em Vertex $a$ is adjacent to $n_3'$, and $b$ has a neighbor in the set $\{c,n_3'\}$.}
\end{equation}
Suppose first that $x$ is adjacent to $y$. 
If $a$ is non-adjacent to $n_3'$, then $y$-$x$-$a$-$n_2$-$n_3'$-$c$  is a $P_6$, a contradiction.  
Moreover, if $b$ is anticomplete to $\{c,n_3'\}$, then $x$-$y$-$b$-$n_2$-$n_3'$-$c$ is a $P_6$, a contradiction. 
So we may assume that $x$ is non-adjacent to $y$, and thus, by the choice of $x$ and $y$, deduce that $X \setminus N(b)$ is anticomplete to $Y \setminus N(a)$. 

Now it follows from \eqref{eqn:XYinT(4)} that every $z \in  V(G)  \setminus (A_1 \cup A_2 \cup A_3 \cup X \cup Y)$ that is anticomplete to $A_1$ and that has a neighbor in $(X \setminus N(b)) \cup (Y \setminus N(a))$ is already complete to $(X \setminus N(b)) \cup (Y \setminus N(a))$. 
By assumption \ref{thm:XYinTasp2}.(ii), we deduce that $X \setminus N(b) = \{x\}$, and similarly $Y\setminus N(a)=\{y\}$. 
Moreover, since $x$ is non-adjacent to $y$, it follows from 
assumption \ref{thm:XYinTasp4} and \eqref{eqn:XYinT(4)}  that 
there exist $x' \in X \cap N(b)$ and $y' \in Y \cap N(a)$ such that $xy'$ and $yx'$ are edges.
Now if $a$ is non-adjacent to $n_3'$, then $y$-$x'$-$a$-$n_2$-$n_3'$-$c$ is a $P_6$, and if $b$ is anticomplete to $\{c,n_3'\}$, then $x$-$y'$-$b$-$n_2$-$n_3'$-$c$ is a $P_6$, in both cases a contradiction. 
This proves \eqref{eqn:XYinT(5)}.

If $b$ is adjacent to $n_3'$, then \eqref{thm:XYinToutcome2} holds, and thus we may assume the opposite.
By \eqref{eqn:XYinT(5)}, $b$ is adjacent to $c$. 
Since $x$-$a$-$n_3'$-$c$-$b$-$y$ is not a $P_6$, we deduce that $x$ is adjacent to $y$. 
Similarly, $X \setminus N(b)$  is complete to $Y \setminus N(a)$.

Let $d=n_1(n_3')$. 
Then $t(d) \leq t(n_2) < t(a)$, and therefore $a \neq d$.  
Since $d$-$n_3'$-$a$-$x$-$y$-$b$ is not a $P_6$, we deduce that $d$ is complete to one of $X \setminus N(a)$  and $Y \setminus N(b)$.

By \eqref{eqn:XYinT(2)}, $d$ is not complete to both $X \setminus N(b)$ and $Y \setminus N(a)$.
Suppose first that $d$ is complete to  $X \setminus N(b)$. 
Then there is some $y' \in Y \setminus N(a)$ that is non-adjacent to $d$. 
Since $n_3'$-$d$-$x$-$y'$-$b$-$n_3$ is not a $P_6$, we deduce that $d$ is adjacent to $n_3$. 
Since $t(d) < t(a)$, $d$ is not complete to $X$, and so there is $x' \in X \cap N(b)$ that is non-adjacent to $d$.
Since $x'$-$b$-$c$-$n_3'$-$d$-$x$ is not a $P_6$, $d$ is adjacent to $c$. 
But $dx,bx'$ is a 2-edge matching between $\{d,b\}$ and $X$, and $\{d,b\}$ is complete to $\{c,n_3\}$, contrary to Claim~\ref{XinT}\eqref{XinT(a)}. 

This proves that $d$ is not complete to $X \setminus N(b)$, and thus $d$ is complete to $Y \setminus N(a)$ and has a non-neighbor $x' \in X \setminus N(b)$. 
Suppose that $d$ is non-adjacent to $n_2$. 
Since $t(n_2(d)) < t(d) \leq  t(n_2)$, we  deduce that $t(n_2(d)) < t(n_2)$, and $a$ is non-adjacent to $n_2(d)$ (since $n_2=n_2(a)$). 
But now $n_2(d)$-$d$-$y$-$x'$-$a$-$n_2$ is a $P_6$, a contradiction. 
This proves that $d$ is adjacent to $n_2$.

Since $\{a,d\}$ is complete to $\{n_2,n_3'\}$, we deduce that there is no 2-edge matching between $Y$ and $\{a,d\}$, by Claim~\ref{XinT}\eqref{XinT(a)}. 
But then $d$ is complete to $Y$, and \eqref{thm:XYinToutcome2} holds, since $n_2$ is adjacent to $n_3'$. 
This completes the proof.
\end{proof}

\begin{claim}\label{thm:G'}
Let $G'$ be the graph obtained from $G$ by contracting $(A_1,A_2,A_3)$ to a triangle $a_1a_2a_3$. 
Let $H'$ be an induced subgraph of $G'$ with $a_1 \in V(H')$. 
Assume that no two non-adjacent neighbors of $a_1$ dominate each other in $H'$. 
Moreover, assume also that for every $v \in V(H')$, either 
\begin{enumerate}
	\item\label{thm:G'asp1} $N_{H'}(v)=X'\cup Y'$, each of $X'$, $Y'$ is stable,
	\begin{enumerate}[(i)]
		\item\label{thm:G'asp1i} for every $x \in X'$ and $y \in Y'$, either
		\begin{enumerate}[(A)]
			\item $x$ is adjacent to $y$,
			\item $x$ has a neighbor in $V(H') \setminus (N_{H'}(v)\cup \{v\})$, or
			\item $y$ has a neighbor in $V(H') \setminus (N_{H'}(v)\cup \{v\})$;
		\end{enumerate}
		\item\label{thm:G'asp1ii} for every $x,x' \in X$ there exists  $p  \in V(H') \setminus \{v\}$ such that $p$ is non-adjacent to $v$, and $p$ is adjacent to exactly one of $x,x'$;	
		\item\label{thm:G'asp1iii} \eqref{thm:G'asp1ii} holds for $Y$ in an analogous way.
	\end{enumerate}
	\item\label{thm:G'asp2} $N_{H'}(v)$ is a triangle, or
	\item\label{thm:G'asp3} $N_{H'}(v)$ induces a $C_5$.
\end{enumerate}
Then either 
\begin{enumerate}[(a)]
	\item\label{thm:G'outcome6.1} some $a \in A_1$ is complete to $N_H'(a_1) \setminus \{a_2,a_3\}$; or
	\item\label{thm:G'outcome6.2} assumption \ref{thm:G'asp1} holds, and no vertex of $A_1$ is complete to $N_{H'}(a_1) \setminus \{a_2,a_3\}$, and there exist $a,b \in A_1$, $n_2 \in A_2$, and $n_3 \in A_3$ such that $a$ is complete to $X' \setminus \{a_2,a_3\}$, $b$ is complete to $Y' \setminus \{a_2,a_3\}$, $\{a,b\}$ is complete to $\{n_2,n_3\}$, and $n_2$ is adjacent to $n_3$;
	\item\label{thm:G'outcome6.3} assumption~\ref{thm:G'asp2} or \ref{thm:G'asp3} holds, and $G$ contains a non-3-colorable graph with seven or eight vertices; or
	\item\label{thm:G'outcome6.4} assumption~\ref{thm:G'asp2} or \ref{thm:G'asp3} holds, there exists a set $A \subseteq A_1$, with $|A| \leq 3$, $n_2 \in A_2$, and $n_3 \in A_3$ such that every vertex of $N_{H'}(a_1)$ has a neighbor in $A$, $A$ is complete to $\{n_2,n_3\}$, and $n_2$ is adjacent to $n_3$.
\end{enumerate}
Moreover, suppose $a_2,a_3$ are in $V(H')$. 
If \eqref{thm:G'outcome6.1} holds, let $H=G|((V(H') \setminus \{a_1\})\cup \{a\})$. Then $H$ is isomorphic to $H'$.
If \eqref{thm:G'outcome6.2} holds, let  $H=G|((V(H') \setminus \{a_1\}) \cup \{a,b,n_1,n_2\})$. 
Then in every coloring of $H$, $a$ and $b$ have  the same color.
If \eqref{thm:G'outcome6.4} holds, let  $H=G|((V(H') \setminus \{a_1\}) \cup A \cup \{n_1,n_2\})$. 
Then in every coloring of $H$, $A$ is monochromatic.

In all cases, $H$ is 3-colorable if and only if $H'$ is.
\end{claim}
\begin{proof} 
Suppose \eqref{thm:G'outcome6.1} does not hold. 

Assume first that assumption \ref{thm:G'asp1} holds for $a_1$.  
Let $X=X' \setminus \{a_2,a_3\}$ and $Y=Y'\setminus \{a_2,a_3\}$. 
We now quickly check that the assumptions of Claim~\ref{thm:XYinT} hold for $A_1,X,Y$ (in $G$).
\begin{itemize}
	\item Every vertex $v \in X\cup Y$ has a neighbor in $A_1$, since every such $v$ is adjacent to $a_1$ in $H'$.
	\item Assumption \ref{thm:XYinTasp1} of Claim~\ref{thm:XYinT} follows from assumption \ref{thm:G'asp1}.(i) of Claim~\ref{thm:G'}.
	\item Assumption \ref{thm:XYinTasp2} holds since there is such a $p$ by assumption \ref{thm:G'asp1}.(ii) of Claim~\ref{thm:G'}.
	Since $p$ is non-adjacent to $a_1$, we deduce that $p \not \in \{a_2,a_3\}$, and so $p \in V(G) \setminus (A_1\cup A_2\cup A_3)$, as desired.
	\item Assumption \ref{thm:XYinTasp3} of Claim~\ref{thm:XYinT} follows analogously.
	\item Assumption \ref{thm:XYinTasp4} of Claim~\ref{thm:XYinT} is seen like this: $N(u)\setminus \{a_1\}$ and $N(v)\setminus \{a_1\}$ are incomparable in $H'$, and $\{u,v\}$ is anticomplete to $\{a_2,a_3\}$ by the maximality of the tripod.
\end{itemize}
Now Claim~\ref{thm:G'} follows from Claim~\ref{thm:XYinT}. 

Next assume that assumption \ref{thm:G'asp2} holds for $a_1$, and $N(a_1)=\{x_1,x_2,x_3\}$. We claim that $\{x_1,x_2,x_3\} \cap \{a_2,a_3\} = \emptyset$. 
Suppose not; we may assume that $x_1=a_2$, and $x_3 \not \in \{a_2,a_3\}$.
Then, in $G$, $x_3$ has both a neighbor in $A_1$ and a neighbor in $A_2$, 
contrary to the maximality of the tripod.

Assume first that there exist $b_1,b_2,b_3 \in A_1$ such that $b_i$ is complete to $\{x_j,x_k\}$ (where $\{1,2,3\}=\{i,j,k\}$). 
Since \eqref{thm:G'outcome6.1} does not hold, $b_i$ is non-adjacent to $x_i$, $i=1,2,3$.  
If some $n_2 \in A_2$ is complete to $\{b_1,b_2,b_3\}$, then \eqref{thm:G'outcome6.3} holds.
So we may assume that there is a 2-edge matching from $A_2$ to $\{b_1,b_2,b_3\}$, say $n_2b_1,n_2'b_2$.
But then $n_2$-$b_1$-$x_2$-$x_1$-$b_2$-$n_2'$ is a $P_6$, a contradiction. 
So we may assume that no vertex of $A_1$ is adjacent to both $x_1$ and $x_2$. 
For $i=1,2$, let $c_i'$ be the smallest vertex in $A_1$ adjacent to $x_i$ w.r.t.~their value of $t$. 
By Claim~\ref{thm:XYinT} applied with $X=\{x_1\}$ and $Y=\{x_2\}$, and since no vertex of $A_1$ is adjacent to both $x_1$ and $x_2$, we deduce that there exist a neighbor $c_i$ of $x_i$, and vertices $n_2 \in A_2$ and $n_3 \in A_3$, such that $\{c_1,c_2\}$ is complete to $\{n_2,n_3\}$, and $n_2$ is adjacent to $n_3$.
If $x_3$ is adjacent to one of $c_1,c_2$, then \eqref{thm:G'outcome6.4} holds, so we may suppose this is not the case. 
Let $c_3$ be a neighbor of $x_3$ in $A$. 
We may assume that $c_3$ is non-adjacent to $x_1$.  
Now $c_3$-$x_3$-$x_1$-$c_1$-$n_2$-$c_2$ is not a $P_6$, and so $c_3$ is adjacent to $n_2$. 
Similarly, $c_3$ is adjacent to $n_3$. 
But now \eqref{thm:G'outcome6.4} holds. 
This finishes the case when assumption \ref{thm:G'asp2} holds.

Finally, assume that \ref{thm:G'asp3} holds. 
Let $N_{H'}(a_1)=\{x_1,\ldots,x_5\}=X$, where $x_1$-$x_2$-$\ldots$-$x_5$-$x_1$ is a $C_5$. 
Since $H'|X$ is connected, the maximality of the tripod implies that 
$\{a_2,a_3\} \cap X = \emptyset$. 
Let $A$ be a minimum size subset of $A_1$ such that each of $x_1,\ldots,x_5$ has a neighbor in $A$. 
Since every $a \in A$ has a neighbor in $A_2$, we deduce that every $a \in A$ has two non-adjacent neighbors in $X$, due to $P_6$-freeness. 
We may assume that $|A|>1$, or \eqref{thm:G'outcome6.3} holds, and so every $a \in A$ is either a \emph{clone} (i.e., has two non-adjacent or three consecutive neighbors in $X$), a \emph{star} (i.e., has four neighbors in $X$), or a \emph{pyramid} for $G|X$ (i.e., has three neighbors in $X$, one of which is non-adjacent to the other two).

Suppose some $a \in A$ is a clone. 
We may assume $a$ is adjacent to $x_2$ and $x_5$. 
If $a$ is mixed on $A_2\cup A_3$, then, since $T_1$ is connected, there is an induced path $a$-$p$-$q$ where $p,q \in A_2\cup A_3$.
There is also an induced path $a$-$x_2$-$x_3$-$x_4$, so $q$-$p$-$a$-$x_2$-$x_3$-$x_4$ is a $P_6$, a contradiction. 
So $a$ is complete to $A_2\cup A_3$. 
If at most one vertex of $A$ is not a clone and $|A| \leq 3$, then 
(using $n_2,n_3$ from  Claim~\ref{thm:2.5} applied to the unique vertex of $A$ that is not a clone, if one exist), outcome \eqref{thm:G'outcome6.4} holds.
So we may assume that if $|A| \leq 3$, then there are at least two non-clones in $A$. 

We claim that $a$ is adjacent to $x_1$. 
Suppose that this is false, and let $b \in A$ be adjacent to $x_1$. 
By the minimality of $A$, $b$ is not complete to $\{x_2,x_5\}$. 
Since $b$ has two non-adjacent neighbors in $X$, by symmetry we may assume that $b$ is adjacent to $x_4$. 
If $b$ is  adjacent to $x_3$, then, by the minimality of $A$, $A=\{a,b\}$ and $b$ is the unique non-clone in $A$, so $b$ is non-adjacent to $x_3$. 
Now $|A\setminus \{a,b\}|=1$, and so $b$ is not a clone. 
Therefore $b$ is adjacent to $x_2$. 

By the minimality of $A$, $b$ is  non-adjacent to $x_5$. 
Let $c \in A$ be adjacent to $x_3$. Then $A=\{a,b,c\}$. By 
the minimality of $A$, $c$ is non-adjacent to $x_5$, and to at least one of $x_1,x_4$. But 
now $c$ is a clone, and $b$ is the unique non-clone in $A$, a contradiction. So $a$ 
is adjacent to $x_1$. This implies that $A=\{a,b,c\}$, $b$ is adjacent to $x_4$ but not to $x_3$,
$c$ is adjacent to $x_3$ but not $x_4$, neither of $b,c$ is a clone, and 
(by the minimality of $A$) no vertex
of $A_1$ is complete to $\{x_3,x_4\}$. By Claim~\ref{thm:XYinT}, there exist $b',c' \in A_1$,
$n_2 \in A_2$ and $n_3 \in A_3$, such that $b'x_4$ and $c'x_5$ are edges, $n_2$
is adjacent to $n_3$, and $\{b',c'\}$ is complete to $\{n_2,n_3\}$. Now \eqref{thm:G'outcome6.4} holds.
So we may assume that $A$ does not contain a clone.
\begin{equation}\label{eqn:Thm6eqn(2)}
\begin{minipage}[c]{0.8\textwidth}\em
If $A=\{a,b\}$ and there exist $x,y,z \in X$ such that $z$-$a$-$x$-$y$-$b$ or $a$-$x$-$y$-$b$-$z$ is an induced path, then \eqref{thm:G'outcome6.2} holds.
\end{minipage}\ignorespacesafterend 
\end{equation}
Since $p$-$a$-$x$-$y$-$b$-$q$ is not a $P_6$ for any $p,q \in A_2$, we deduce
that either $N(a) \cap A_2 \subseteq N(b) \cap A_2$, or 
$N(b) \cap A_2 \subseteq N(a) \cap A_2$, and the same holds in $A_3$. Since we may assume
\eqref{thm:G'outcome6.2} does not hold, Claim~\ref{thm:2.5} implies that, up to symmetry, there exist 
$n_2 \in N(a) \cap A_2$ and $n_3 \in N(b) \cap A_3$ such that $a$ is non-adjacent to $n_3$,
and $b$ is non-adjacent to $a_2$. Then $n_2$ is adjacent to $n_3$ (or $n_2$-$a$-$x$-$y$-$b$-$n_3$ is
a $P_6$). But now $z$-$a$-$n_2$-$n_3$-$b$-$y$ or $z$-$b$-$n_3$-$n_2$-$a$-$x$ is a $P_6$, a contradiction. This 
proves \eqref{eqn:Thm6eqn(2)}.

Suppose some $a \in A$ is a star, say $a$ is adjacent to $x_1,\ldots,x_4$, and not to $x_5$.
Let $b \in A$ be adjacent to $x_5$. Then we know that $A=\{a,b\}$. If $b$ is adjacent to both $x_1$ and 
$x_4$, then \eqref{thm:G'outcome6.3} holds, and so we may assume that $b$ is non-adjacent to $x_1$.  Since $b$ is not
a clone, $b$ is adjacent to $x_2$. If $b$ is adjacent to $x_3$, then \eqref{thm:G'outcome6.3} holds, so
$b$ is non-adjacent to $x_3$; since $b$ is not a clone, $b$ is adjacent to $x_4$. But
now \eqref{eqn:Thm6eqn(2)} holds with $x=x_1$, $y=x_5$ and $z=x_3$. So we may assume that no $a \in A$
is a star, and so every vertex of $A$ is a pyramid.

Let $a \in A$.
We may assume that $a$ is adjacent to $x_1,x_3,x_4$ and not to $x_2,x_5$. 
Let $b \in A$ be adjacent to $x_2$. 
If $N(b) \cap X=\{x_2,x_4,x_5\}$, then \eqref{thm:G'outcome6.2} holds by \eqref{eqn:Thm6eqn(2)}
applied with $x=x_3$, $y=x_2$ and $z=x_5$. 
If $N(b) \cap X = \{x_2,x_3,x_5\}$, then
we obtain the previous case by exchanging the roles of $a$ and $b$. 
So we may assume that $N(b) \cap X =\{x_1,x_2,x_4\}$.

Hence, there exists $c \in A\setminus \{a,b\}$ adjacent to $x_5$ with $N(c) \cap X=\{x_1,x_3,x_5 \}$.
But now every $x \in X$ has a neighbor in $A\setminus \{a\}$, contrary to the minimality of $A$.
This shows how the statement of Claim~\ref{thm:G'} follows from assumption \ref{thm:G'asp3}, completing the proof.
\end{proof}

\begin{claim}\label{thm:7}
Every graph $H' \in \mathcal L$ satisfies the assumptions of Claim~\ref{thm:G'}. 
\end{claim}
\begin{proof}
Since $H'$ is a minimal obstruction to 3-coloring, $H'$ has no 
dominated vertex, meaning any two neighborhoods of vertices are incomparable. 
Let $v \in V(H')$. 
If $N(v)$ is not bipartite, then $v$ contains a triangle  or $C_5$, and so $V(H')=\{v\}\cup N(v)$, and 
assumptions \ref{thm:G'asp2} or \ref{thm:G'asp3} of Claim~\ref{thm:G'} hold. So $N(v)$ is bipartite with a bipartition 
$(X,Y)$. 

We implemented a straightforward program which we used to verify that assumption~\ref{thm:G'asp1} of Claim~\ref{thm:G'} indeed holds for all 24 $4$-critical $P_6$-free graphs from~Theorem~\ref{thm:N3P6} where $N(v)$ is bipartite. The source code of this program can be downloaded from~\cite{tripodgenerator-site}.
\end{proof}

\begin{claim}\label{thm:8}
Let $G'$ be obtained from G by contracting $(A_1,A_2,A_3)$ to a triangle $a_1a_2a_3$. 
Let $H'$ be an induced subgraph of $G'$, with $a_1,a_2 \in V(H')$. For $i=1,2$,
let $Z_i=N(a_i) \setminus \{a_1,a_2,a_3\}$. 

Assume that
\begin{enumerate}
	\item\label{thm:8A2} no two non-adjacent neighbors of $a_1$ dominate each other, and no two non-adjacent neighbors of $a_2$ dominate each other, and
	\item\label{thm:8A3} $H'|N(a_1)$ and $H'|N(a_2)$ are bipartite.
\end{enumerate}

Then Claim~\ref{thm:G'}\eqref{thm:G'outcome6.1} or Claim~\ref{thm:G'}\eqref{thm:G'outcome6.2} holds for each of $a_1,a_2$.
If Claim~\ref{thm:G'}\eqref{thm:G'outcome6.1} holds for $a_1$, let $c_1$ be the vertex $a$ of Claim~\ref{thm:G'}\eqref{thm:G'outcome6.1}, set $A=\{c_1\}$ and $Z=\emptyset$.
If Claim~\ref{thm:G'}\eqref{thm:G'outcome6.2} holds for $a_1$, let $a,b,n_2(a_1),n_3(a_1)$ be the vertices as in Claim~\ref{thm:G'}\eqref{thm:G'outcome6.2}. 
Moreover, set $A=\{a,b\}$, and $Z= \{n_2(a_1),n_3(a_1)\}$.

If Claim~\ref{thm:G'}\eqref{thm:G'outcome6.1} holds for $a_2$, let $c_2$ be the vertex $a$ of Claim~\ref{thm:G'}\eqref{thm:G'outcome6.1}, set $C=\{c_2\}$, and $W=\emptyset$. 
If Claim~\ref{thm:G'}\eqref{thm:G'outcome6.2} holds for $a_2$, let $c,d,n_1(a_2),n_3(a_2)$ be the vertices as in Claim~\ref{thm:G'}\eqref{thm:G'outcome6.2}, set
$C=\{c,d\}$, and $W=\{n_1(a_2),n_3(a_2)\}$.
$C$.

One of the following holds.
\begin{enumerate}[(a)]
	\item\label{thm:8outcome8.0} Outcome Claim~\ref{thm:G'}\eqref{thm:G'outcome6.1} holds for $a_1$, there is $c \in C$, and an induced path $c_1$-$c'$-$a'$-$c$
in $T_3(t)$ where $t=\max(t(c_1),t(c))$, such that $a'$ is complete to $Z_1$. Or the analogous statement holds for $a_2$.
	\item\label{thm:8outcome8.1} There is an edge between $A$ and $C$.
In this case let $H=(H' \setminus \{a_1,a_2\}) \cup A \cup C \cup Z \cup W$.	
	\item\label{thm:8outcome8.2} In $H'$, there is an induced path $a_1$-$q_1$-$q_2$-$a_2$, and a vertex complete to $\{a_1,q_1,q_2\}$ or to $\{a_2,q_2,q_1\}$.In this case let $H=(H' \setminus \{a_1,a_2\}) \cup A \cup C \cup Z \cup W$.
	\item\label{thm:8outcome8.3} There are adjacent vertices $n_1 \in A_1$ and $n_2 \in  A_2$, such that  
$n_1$ is complete to $C$, $n_2$ is complete to $A$, and some vertex $s \in A_3$ is complete to $A \cup C \cup \{n_1,n_2\}$.
In this case let 
$H=(H' \setminus \{a_1,a_2\}) \cup A \cup C \cup \{n_1,n_2,s\}$.
\end{enumerate}

In each  the cases \eqref{thm:8outcome8.1}, \eqref{thm:8outcome8.2}, \eqref{thm:8outcome8.3}, in every 3-coloring of $H$, $A$ and $C$ are monochromatic, and no
color appears in both $A$ and $C$. 

In all cases, $H$ is 3-colorable if and only if $H'$ is.

\end{claim}
\begin{proof}
By assumption \ref{thm:8A3} Claim~\ref{thm:G'}\eqref{thm:G'outcome6.1} or 
Claim~\ref{thm:G'}\eqref{thm:G'outcome6.2} holds for each of $a_1,a_2$.
We may assume that no vertex of $V(G) \setminus (Z_1\cup A_2\cup A_3)$ has a neighbor in $A_1$, and no vertex
of $V(G) \setminus (Z_2\cup A_1\cup A_3)$ has a neighbor in $A_2$ (otherwise we may delete such vertices from $G$
without changing the hypotheses or the outcomes).

Moreover, we may assume that $A$ is anticomplete to $C$, as otherwise \eqref{thm:8outcome8.1} holds.  
Pick $a \in A$ and $c \in C$.
Let $t=\max(t(a),t(c))$, and let $c$-$a'$-$c'$-$a$ be an induced path from $a$ to $c$ in $T_3(t)$.  
If possible, we choose $a'$ to be complete 
to $C$, and $c'$ complete to $A$. 

Assuming \eqref{thm:8outcome8.0} does not hold, we derive the following.
\begin{equation}\label{eqn:thm8(1)}
\mbox{\em Vertex $a'$ is not complete to  $Z_1$, and $c'$ is not complete to $Z_2$.}
\end{equation}
We also make use of the following fact.
\begin{equation}\label{eqn:thm8(2)}
\mbox{\em Vertex $c'$ is complete to $A$, and $a'$ to $C$.}
\end{equation}
To see this, suppose $c'$ is not complete to $A$. Then $A=\{a,b\}$, and $c'$ is non-adjacent to $b$.
By the choice of $c'$, we deduce that $n_2(a_1)$ is non-adjacent to $a'$ 
(otherwise we may replace $c'$ with $n_2(a_1)$). Now $b$-$n_2(a_1)$-$a$-$c'$-$a'$-$c$ is a $P_6$, a 
contradiction. Similarly, $a'$ is complete to $C$. This proves \eqref{eqn:thm8(2)}.
\begin{equation}\label{eqn:thm8(3)}
\begin{minipage}[c]{0.8\textwidth}\em
Let $p \in Z_1$ be non-adjacent to $a'$. Then $p$ has no neighbor in  
$V(H') \setminus (\{a_1,a_2,a_3\}\cup Z_1\cup Z_2)$, and $p$ has a neighbor $q \in Z_1$.
\end{minipage}\ignorespacesafterend 
\end{equation}
Since $a_2$ does not dominate $p$, $p$ has a neighbor $q \in H'$ non-adjacent to 
$a_2$. Then in $G$, $q$ is anticomplete to $A_2$.  Let $z \in A$ be adjacent to $p$. If $q$ is 
not in $Z_1$, then $q$ is anticomplete to $A_1$, and so, by \eqref{eqn:thm8(2)},  $q$-$p$-$z$-$c'$-$a'$-$c$ is 
a $P_6$ in $G$, a contradiction. This proves \eqref{eqn:thm8(3)}.

By \eqref{eqn:thm8(1)}, \eqref{eqn:thm8(3)} and the symmetry between $A_1$ and $A_2$, there exist $p,q \in Z_1$
and $s,t \in Z_2$ such that $pq,st$ are edges, $a'$ is non-adjacent to $p$, and $c'$
is non-adjacent to $s$. Let $r \in A$ be adjacent to $p$, and let $u \in C$ be
adjacent to $s$. Since $p$-$r$-$c'$-$a'$-$u$-$s$ is not a $P_6$, we deduce that $p$ is adjacent to $s$.

Let $D$ be the following $C_6$: $r$-$c'$-$a'$-$u$-$s$-$p$-$r$.
\begin{equation}\label{eqn:thm8(4)}
\mbox{\em Vertex $p$ is complete to $A$, and $s$ is complete to $C$.}
\end{equation}

Suppose $p$ has a non-neighbor $r' \in A$. Then, since $A$ is anticomplete to $C$, $r'$ is a 
leaf for $D$, in contradiction to Claim~\ref{clm:noleaforhat}. Similarly, $s$ is complete to $C$. This proves \eqref{eqn:thm8(4)}.

By \eqref{eqn:thm8(4)}, we may assume that $r$ is adjacent to $q$, and $u$ is adjacent to $t$.
If $q$ is adjacent to $s$, then \eqref{thm:8outcome8.2} holds, which we may assume not to be the case. 
Similarly, $t$ is non-adjacent to $p$.
Since $q,t$ are not hats for $D$, by Claim~\ref{clm:noleaforhat}, we deduce that $q$ is adjacent  to $a'$, and $t$ to $c'$.

Let $d \in A_3$ be adjacent to $a'$.  If $d$ is non-adjacent to $c'$, then $d$-$a'$-$c'$-$t$-$s$-$p$ is a $P_6$, a contradiction.  So $d$ is adjacent to $c'$.
Since by Claim~\ref{clm:noleaforhat} $d$ is not a hat for $D$, we deduce that $d$ is adjacent to at least one of $r,u$. Suppose that $d$ is adjacent to $r$ and not to $u$. Then 
$p$-$r$-$d$-$a'$-$u$-$s$-$p$ is a $C_6$, and $t$ is a hat for it, again contrary to Claim~\ref{clm:noleaforhat}. This proves that $d$ is complete to $\{r,u\}$. Similarly
$d$ is complete to $A \cup C$ and \eqref{thm:8outcome8.3} holds. 

It is now easy to verify that the last assertion of Claim~\ref{thm:8} holds.
This completes the proof.
\end{proof}

By $W_5$ we denote the graph that is $C_5$ plus a vertex adjacent to all vertices of that $C_5$.

\begin{claim}\label{thm:9} 
Every $H \in \mathcal L$ except $K_4$ and $W_5$
satisfies the assumptions of Claim~\ref{thm:8}.
\end{claim}
\begin{proof} 
Let $H \in \mathcal L$. Since $H$ is minimal non-3-colorable, $H$ has no 
dominated vertices, and so assumption \ref{thm:8A2} of Claim~\ref{thm:8} holds. If $H|N(v)$ is not bipartite for some $v \in 
V(H)$, then  $H|N(v)$ contains a triangle or a $C_5$, and so $H=K_4$ or $H=W_5$.
\end{proof}

We can now prove our main statement of this section.

\begin{proof}[Proof of Lemma~\ref{thm:10}]
We may assume that at least one of $a_1,a_2,a_3$ is in $V(H')$. If $|V(H') \cap \{a_1,a_2,a_3\}|=1$,
we are done using Claim~\ref{thm:G'} and Claim~\ref{thm:7}, so we may assume that $|V(H') \cap \{a_1,a_2,a_3\}| \geq 2$.
Note that if  $H'=K_4$, every edge is in a triangle, and if $H'=W_5$, then
every triangle is in a diamond.
Hence, the maximality of $(A_1,A_2,A_3)$ implies that
$H' \neq K_4, W_5$.

By Claim~\ref{thm:9}, $H'$ satisfies the assumptions of Claim~\ref{thm:8}. 
We define the sets $C_i,W_i,N_i,Z_i$ for $i \in \{1,2,3\}$. If $a_i \not \in V(H')$, let $C_i=W_i= Z_i= \emptyset$, and let 
$N_j=\emptyset$ for every $j \neq i$. Now suppose that $a_1 \in V(H')$. 
Let $C_1=A$ be as in 
Claim~\ref{thm:G'}\eqref{thm:G'outcome6.1} or Claim~\ref{thm:G'}\eqref{thm:G'outcome6.2}. (Observe that there may be several possible choices of $C_1$.)
If Claim~\ref{thm:G'}\eqref{thm:G'outcome6.2} holds for
$a_1$, let $W_1=\{n_2,n_3\}$ (in the notation of Claim~\ref{thm:G'}\eqref{thm:G'outcome6.2}), if  If Claim~\ref{thm:G'}\eqref{thm:G'outcome6.1} holds for $a_1$, let $W_1=\emptyset$.
Let $Z_1$ be the set of neighbors of $a_1$ in $H' \setminus \{a_1,a_2,a_3\}$.  
Suppose $a_2 \in V(H')$.
If Claim~\ref{thm:8}\eqref{thm:8outcome8.3} holds for $a_1a_2$, define $N_3=\{n_1,n_2,s\}$ in the notation of Claim~\ref{thm:8}\eqref{thm:8outcome8.3}. If Claim~\ref{thm:8}\eqref{thm:8outcome8.0}, Claim~\ref{thm:8}\eqref{thm:8outcome8.1} or
Claim~\ref{thm:8}\eqref{thm:8outcome8.2} holds for $a_1a_2$, let
$N_3=\emptyset$. Note that in all  cases $|N_3|\leq 3$.
Define $C_2,C_3,W_2,W_3,N_1,N_2,Z_2,Z_3$ similarly.

Next we define the sets $D_1,D_2,D_3$. If Claim~\ref{thm:8}\eqref{thm:8outcome8.3} holds for the pair $a_1a_2$ or for the pair $a_1a_3$,  we set 
$D_1=N_2 \cup N_3$, and otherwise we set $D_1=W_1$.

As usual, we may assume $V(G)=A_1\cup A_2\cup A_3\cup (V(H') \setminus \{a_1,a_2,a_3\})$.

We analyze the possible outcomes of Claim~\ref{thm:8}.

Let us call outcomes
Claim~\ref{thm:8}\eqref{thm:8outcome8.1}, Claim~\ref{thm:8}\eqref{thm:8outcome8.2}, Claim~\ref{thm:8}\eqref{thm:8outcome8.3}  {\em good}. Suppose first that a good outcome holds for each pair $a_1a_2,a_2a_3,a_1a_3$ contained in $V(H')$. 

Let
$$H=(H' \setminus \{a_1,a_2,a_3\}) \cup (C_1 \cup C_2 \cup C_3 \cup D_1 \cup D_2 \cup D_3).$$

\begin{itemize}
	\item If Claim~\ref{thm:G'}\eqref{thm:G'outcome6.1} holds for each of $a_1,a_2,a_3$, then $|V(H)| \leq |V(H')|+9$, as follows.

We observe that in this case  $|C_i|=1$ and $|W_i|=0$ for each $i$,
and therefore $D_1 \cup D_2 \cup D_3=N_1 \cup N_2 \cup N_3$.
Since $|N_i| \leq 3$ for every $i$, we have that 
$|V(H)| \leq |V(H')|-3+3+9=|V(H')|+9$.

	\item If Claim~\ref{thm:G'}\eqref{thm:G'outcome6.1} holds for exactly 
two of $a_1,a_2,a_3$, then $|V(H)| \leq |V(H')|+12$, as follows.

Assume Claim~\ref{thm:G'}\eqref{thm:G'outcome6.1} holds for $a_1$ and $a_2$.
Then   Claim~\ref{thm:G'}\eqref{thm:G'outcome6.2} holds for $a_3$.
It follows that  $|C_1|=|C_2|=1$, $W_1=W_2=\emptyset$,
and $|C_3|=|W_3|=2$. Therefore 
$D_1 \cup D_2 \cup D_3 \subseteq N_1 \cup N_2 \cup N_3 \cup W_3$, and so
$|D_1 \cup D_2 \cup D_3| \leq 11$. Consequently, 
$|V(H)| \leq |V(H')|-3+4+11=|V(H')|+12$.

	\item If Claim~\ref{thm:G'}\eqref{thm:G'outcome6.1} holds for exactly one of $a_1,a_2,a_3$, then $|V(H)| \leq |V(H')|+12$, as follows.

Assume Claim~\ref{thm:G'}\eqref{thm:G'outcome6.1} holds for $a_1$. 
Then Claim~\ref{thm:G'}\eqref{thm:G'outcome6.2} holds for exactly $a_2,a_3$.
It follows  that  $|C_1|=1$, $W_1=\emptyset$ and  $|C_2|=|W_2|=|C_3|=|W_3|=2$.
We claim that $D_1 \cup D_2 \cup D_3 \leq 10$.
Since $W_1 = \emptyset$, we deduce that $D_1 \subseteq N_2 \cup N_3$. 
If $N_1 \neq \emptyset$. Then $D_2 \cup D_3  \subseteq N_1 \cup N_2 \cup N_3$, and so, since $|N_i| \leq 3$ for every $i$,  
$|D_1 \cup D_2 \cup D_3| \leq |N_1 \cup N_2 \cup N_3| \leq 9$.
Thus we may assume that $N_1=\emptyset$. Then
$D_1 \cup D_2 \cup D_3 \subseteq N_2 \cup N_3 \cup W_2 \cup W_3$,
and again $|D_1 \cup D_2 \cup D_3| \leq 10$.
Consequently,
$|V(H)| \leq |V(H')|-3+5+10=|V(H')|+12$.

	\item If Claim~\ref{thm:G'}\eqref{thm:G'outcome6.2} holds for all of $a_1,a_2,a_3$, then $|V(H)| \leq |V(H')|+12$, as follows.

Since  Claim~\ref{thm:G'}\eqref{thm:G'outcome6.2} holds for each of 
$a_1,a_2,a_3$, it follows that  for every $|C_i|=|W_i|=2$ for every $i$. 

We show that $|D_1 \cup D_2 \cup D_3| \leq 9$. 
Suppose first that $N_1 \neq \emptyset$ and $N_2 \neq \emptyset$. Then
$D_2 \cup D_3 \cup D_3 \subseteq N_1 \cup N_2 \cup N_3$, and the claim follows
since $|N_i| \leq 3$ for every $i$. This we may assume that 
$N_2=N_3 = \emptyset$. 

Next assume that $N_1 \neq \emptyset$. Then $D_2=D_3=N_1$, and
$D_1=W_1$, and so $|D_1 \cup D_2 \cup D_3| \leq 5$.

Finally, if $N_1=N_2=N_3=\emptyset$, then $D_i=W_i$ for every $i$, and
thus $|D_1 \cup D_2 \cup D_3| \leq 6$. 

Thus, in all cases $|V(H)| \leq |V(H')|-3+6+9=|V(H')|+12$.
\end{itemize}

Now we may assume that for 
least one of the pairs $a_1a_2, a_2a_3, a_1a_3$ contained in $H'$ no good outcome holds. Consequently, for at
least one of the pairs $a_1a_2, a_2a_3, a_1a_3$ contained in $H'$
Claim~\ref{thm:8}\eqref{thm:8outcome8.0} holds.
Observe that  this is true  for every choice of $C_1,C_2,C_3$ as above.

Suppose  that  $V(H') \cap \{a_1,a_2,a_3\} =\{a_1,a_2\}$. Then 
Claim~\ref{thm:8}\eqref{thm:8outcome8.0} holds for the pair $a_1a_2$. In
the notation of Claim~\ref{thm:8}, we may assume that $a'$ is complete to 
$N_{H'}(a_1)$. Then replacing $C_1$ with $\{a'\}$, we observe that
outcome Claim~\ref{thm:8}\eqref{thm:8outcome8.1} holds for the pair
$a_1a_2$, a contradiction. Thus we may assume that $a_1,a_2,a_3 \in V(H')$.

\begin{equation}\label{eqn:thm10(1)}
\begin{minipage}[c]{0.8\textwidth}\em
Permuting the indices if necessary, there exist $b_2,b_3 \in A_1$,
and $C_2 \subseteq A_2$, $C_3 \subseteq A_3$ such that the following holds.
\begin{itemize}
	\item $\{b_2,b_3\}$ is complete to $Z_1$,
	\item $C_2$ and $C_3$ are as in Claim~\ref{thm:G'}\eqref{thm:G'outcome6.1} or Claim~\ref{thm:G'}\eqref{thm:G'outcome6.2},
	\item $b_2$ has a neighbor in $C_2$ and none in $C_3$, 
	\item $b_3$ has a neighbor in $C_3$ and none in $C_2$, and
	\item one of the good  outcomes holds for the pair $C_2,C_3$.
	\item $b_2$ and $b_3$ have a common neighbor in $A_2$ or $A_3$.
\end{itemize}
\end{minipage}\ignorespacesafterend
\end{equation}

In order to prove \eqref{eqn:thm10(1)}, we first prove that a certain condition
is sufficient  for \eqref{eqn:thm10(1)}.
\begin{equation}\label{eqn:thm10(1.1)}
\begin{minipage}[c]{0.8\textwidth}\em
If there exist $C_i' \subseteq A_i$ as in Claim~\ref{thm:G'}\eqref{thm:G'outcome6.1} or Claim~\ref{thm:G'}\eqref{thm:G'outcome6.2} such that there is an 
edge between $C_1'$ and $C_2'$, and an edge between $C_2'$ and $C_3'$, then \eqref{eqn:thm10(1)} holds.
\end{minipage}\ignorespacesafterend
\end{equation}
To see this, apply Claim~\ref{thm:8} with $A=C_1'$ and $C=C_3'$.
If one of the good outcomes holds, then a good outcome holds for all three pairs among $C_1',C_2',C_3'$, and so we may assume that this is not the case. 
There is symmetry between $C_1'$ and $C_3'$, so we may assume that $|C_1'|=1$ and that there is an induced path 
$c_1'$-$c_3''$-$c_1''$-$c_3'$ in $T_2$, where $\{c_1'\}=C_1'$, $c_3' \in C_3'$, and $c_1''$ is complete to $Z_1$. 
If $c_1''$ has a neighbor in $C_2$, or $c_1'$ has a neighbor in $C_3$, then a good outcome holds for all pairs among $\{c_1''\},C_2',C_3'$ or $C_1',C_2',C_3'$.
Hence, we may assume that this is not the case. Now \eqref{eqn:thm10(1)} holds, and this proves \eqref{eqn:thm10(1.1)}.

We may assume that Claim~\ref{thm:8}\eqref{thm:8outcome8.0} holds for the pair $C_2,C_3$.
By modifying $C_2,C_3$ we may assume that there is an edge between $C_2$ and $C_3$ and outcome Claim~\ref{thm:8}\eqref{thm:8outcome8.1} holds for $(C_2,C_3)$. 
If a good outcome holds for both $C_1,C_2$, and $C_1,C_3$, then a good outcome holds for all three pairs, so we may assume that this is not the case. 

So, assume that outcome \eqref{thm:8outcome8.0} holds when Claim~\ref{thm:8} is applied to $C_1,C_2$.  
If there is $c_1 \in A_1$ that is complete to $Z_1$ and has a neighbor in $C_2$, then \eqref{eqn:thm10(1)} holds by \eqref{eqn:thm10(1.1)}.
So we may assume that there is a vertex $c_2' \in A_2$ that is complete to $Z_2$, and an induced path
$c_1$-$c_2'$-$c_1'$-$c_2$ in $T_3$, where $c_1 \in C_1$ and $C_2=\{c_2\}$. 
If a good outcome holds for $C_1,C_3$, then either \eqref{eqn:thm10(1)} holds, or a good outcome holds for 
all three pairs among $C_1,\{c_2\},C_3$ or $C_1,\{c_2'\},C_3$.

So, we may assume that Claim~\ref{thm:8}\eqref{thm:8outcome8.0} holds for $C_1,C_3$. 
By the symmetry between $C_1$ and $C_3$, we may assume that there is $d_1 \in A_1$ and an induced path
$c_1$-$c_3'$-$d_1$-$c_3$ where $c_3 \in C_3$, $C_1=\{c_1\}$, and $d_1$ is complete to $Z_1$.
But now there is an edge between $C_3$ and $\{d_1\}$, and between $C_3$ and $C_2$, and \eqref{eqn:thm10(1)}
follows from \eqref{eqn:thm10(1.1)}. This proves \eqref{eqn:thm10(1)}.

If Claim~\ref{thm:8}\eqref{thm:8outcome8.1} or Claim~\ref{thm:8}\eqref{thm:8outcome8.2} holds for the pair $C_2,C_3$, let $$H=G|((V(H') \setminus \{a_1,a_2,a_3\}) \cup \{b_2,b_3\}\cup C_2\cup C_3\cup W_2\cup W_3),$$
and let $$H''= G|((V(H') \setminus \{a_1,a_2,a_3\}) \cup  \{b_2,b_3\}\cup C_2\cup C_3).$$
If Claim~\ref{thm:8}\eqref{thm:8outcome8.3} holds for the pair $C_2,C_3$, let  $$H=G|((V(H') \setminus \{a_1,a_2,a_3\}) \cup \{b_2,b_3\}\cup C_2\cup C_3\cup N_1),$$ and let $$H''= G|((V(H') \setminus \{a_1,a_2,a_3\})\cup \{b_2,b_3\}\cup C_2\cup C_3).$$ 
Then $|V(H)| \leq |V(H')|+7$, and so we may assume that $H$ is 3-colorable.

Let us call a 3-coloring of $H''$ \emph{promising} if $C_2$ is monochromatic, $C_3$ is monochromatic, and no color appears in both of $C_2,C_3$.
We observe that by Claim~\ref{thm:G'} and Claim~\ref{thm:8}, every 3-coloring of $H$ gives a promising 3-coloring 
of $H''$.
Since $H'$ is not 3-colorable, in every promising coloring of $H''$ the vertices $b_2$ and $b_3$ receive different colors.

Let $c$ be a 3-coloring of $H$. 
We may assume that $c(b_i)=i$, $c$ is constantly 1 or 3 on $C_2$, and $c$ is constantly 1 or 2 on $C_3$. 
Then $c(z)=1$ for every $z \in Z_1$.  
If $c$ is 1 on $C_2$, then we recolor $b_2$ with color 3, and get a coloring of $H'$, a contradiction. 
So we may assume that $c$ is $3$ on $C_2$, and $c$ is 2 on $C_3$. 
If no vertex of $Z_2$ has color 1, we recolor $C_2$ with color 1, and recolor $b_2$ with color 3.
We obtain coloring of $H$ with $b_2,b_3$ colored in the same color, a contradiction. 
So, for some $z_2 \in Z_2$, $c(z_2)=1$. 
Similarly, for some $z_3 \in Z_3$, $c(z_3)=1$. 

For $i=2,3$ let $Z_i'$ be the set of all vertices $z \in Z_i$ with $c(z_i)=1$.
Then $Z_1\cup Z_2'\cup Z_3'$ is a stable set. 
Let $c_i \in C_i$ be adjacent to $b_i$.

\begin{equation}\label{eqn:thm10(2)}
\mbox{\em $Z_2'$ is anticomplete to $V(G) \setminus (Z_2\cup A_2)$.}
\end{equation}
Suppose $p \in V(G) \setminus (Z_2\cup A_2)$ has a neighbor $z_2 \in Z_2'$. 
Then $p \not \in Z_1$. 
Let $c_2' \in C_2$ be adjacent to $z_2$. 
Suppose first that $b_2$ is non-adjacent to $c_2'$.
Then $c_2' \neq c_2$. 
Let $n_1 \in A_1$ be complete to $\{c_2,c_2'\}$, a possible choice by \ref{thm:G'outcome6.2}.
Now $p$-$z_2$-$c_2'$-$n_1$-$c_2$-$b_2$ is a $P_6$, a contradiction. 
So $c_2'$ is adjacent to $b_2$.
Let $n_2 \in A_2$ be adjacent to $b_2$ and $b_3$ (as in \eqref{eqn:thm10(1)}, with the roles of $A_2$ and $A_3$ exchanged).  
Then $p$-$z_2$-$c_2'$-$b_2$-$n_2$-$b_3$ is a $P_6$, again a contradiction. 
This proves \eqref{eqn:thm10(2)}.

Now, by \eqref{eqn:thm10(2)}, we can recolor $H''$ by putting $c'(C_2)=1$ and $c'(Z_2')=3$, $c'(b_2)=3$,
which yields a 3-coloring of $H'$, a contradiction. 
This completes the proof.
\end{proof}

In Section~\ref{sec:proof-of-main-result} we use Lemma~\ref{thm:10} (together with Lemma~\ref{thm:uncontractingK4}, which is the analogue of Lemma~\ref{thm:10} for the case when $H'=K_4$)  to prove the main result of the paper.

\section{Obstructions that are 1-vertex extensions of a tripod}\label{sec:tripod-generator}

In this section, we prove the following statement.

\begin{lemma}\label{thm:uncontractingK4}
Let $G$ be a 4-critical $P_6$-free graph.
Assume that there is a tripod $T=(A_1,A_2,A_3)$ in $G$ and some vertex $x$ which has a neighbor in each $A_i$, $i=1,2,3$.
Then $|V(G)| \le 18$.
\end{lemma}

To see this, let $G$, $T=(A_1,A_2,A_3)$, and $x$ be as in Lemma~\ref{thm:uncontractingK4}.´
Let $a_1,a_2,a_3$ be the root of $T$.
It is clear that $V(G)=V(T)\cup \{x\}$.
We call $G$ a \emph{1-vertex extension of a tripod}.

\subsection{Preparation}

We may assume that the ordering $A_1 \cup A_2 \cup A_3 = \{v_1,\ldots,v_k\}$ has the following property.

\begin{claim}\label{thm:special-ordering}
Let $u \in A_\ell$ and $v \in A_k$ for some $\ell,k \in \{1,2,3\}$.
Moreover, let $\{\ell,\ell',\ell''\} = \{1,2,3\}$ and $\{k,k',k''\} = \{1,2,3\}$.
Assume that $\max(t(n_{k'}(v)),t(n_{k''}(v))) < \max(t(n_{\ell'}(u)),t(n_{\ell''}(u)))$.
Then $t(v) < t(u)$. 
\end{claim}

Let $b_i$ be the neighbor of $x$ in $A_i$ with $t(b_i)$ maximum, for all $i=1,2,3$.
We may assume that $t(b_1) > t(b_2) > t(b_3)$.

\begin{claim}\label{thm:x-has-only-one-neighbor}
We may assume that $N(x)\cap A_1 = \{b_1\}$ and $N(x)\cap A_i = \{b_i\}$ for some $i\in \{2,3\}$.
\end{claim}
\begin{proof}
Since $G|(V(T(t(b_1))) \cup \{x\})$ is 4-chromatic we know that $V(G) = V(T(t(b_1))) \cup \{x\}$.
In particular, $N(x) \cap A_1 = \{b_1\}$.

To see the second statement, assume that $|N(x) \cap A_2|, |N(x) \cap A_3| \ge 2$.
Suppose for a contradiction that $|N(b_1)\cap A_2|, |N(b_1)\cap A_3| \ge 2$, and let $u$ be the vertex in the set $\{b_2,b_3,n_2(b_1),n_3(b_1)\}$ with $t(u)$ maximum.
Then $G-u$ is still 4-chromatic, a contradiction.

So we may assume that $|N(b_1) \cap A_i| = 1$ for some $i \in \{2,3\}$.
Note that $T'=(A_1 \setminus \{b_1\}\cup \{x\},A_2,A_3)$ is a tripod.
Consequently, $b_1$ has neighbors in all three classes of $T'$. 
Since $|N(b_1)\cap (A_1\cup \{x\})|=|N(b_1)\cap A_i|=1$, we are done.
\end{proof}

\subsection{The enumeration algorithm}

Consider the following way of traversing the tripod $T$.
Initially, the vertices $b_1,b_2,b_3$ are labeled \emph{active}, and all other vertices are unlabeled.
Then, we label the vertices $a_1,a_2,a_3$ as \emph{inactive}.
Consequently, if $b_3=a_3$, say, then $b_3$ is labeled inactive.

Iteratively, pick an active vertex, say $u \in A_i$ with $\{i,j,k\}=\{1,2,3\}$.
Make $n_j(u)$ and $n_k(u)$ active, unless they are labeled already, whether active or inactive.
Then label $u$ as inactive and re-iterate, picking another active vertex, if possible.

\begin{claim}\label{thm:active-doesnt-matter}
Regardless of which active vertex is picked in the successive steps, this procedure terminates and, moreover, every vertex of $T$ is visited during this procedure.
\end{claim}
\begin{proof}
Clearly this procedure terminates when there is no active vertex left.
Since every vertex is labeled active at most once, this proves the first assertion.

Assume now the procedure has terminated.
The latter assertion follows from the fact that, if $W$ is the collection of inactive vertices, $G|W$ is already a tripod.
Thus, since $b_1,b_2,b_3 \in W$, $G|(W\cup \{x\})$ is 4-chromatic and so $G|(W\cup \{x\}) = G$, due to the choice of $G$.
\end{proof}

Instead of traversing a given tripod, we use this method to enumerate all possible 4-critical $P_6$-free 1-vertex extensions of a tripod.
The idea is to successively generate the possible subgraphs induced by the labeled vertices only. 
This is done by Algorithm~\ref{algo:init-algo2}.
Starting from all relevant graphs on the vertex set $\{x,b_1,b_2,b_3,a_1,a_2,a_3\}$, we iteratively add new vertices, mimicking the iterative labeling procedure mentioned above.
The following list contains all of these start graphs.

\begin{claim}\label{thm:startgraphs}
We may assume that the graph $G':=G|\{x,b_1,b_2,b_3,a_1,a_2,a_3\}$ has the following properties.
\begin{enumerate}[(a)]
	\item If $b_1=a_1$, then $G=G'$ is $K_4$.
	\item If $b_1 \neq a_1$ and $b_2=a_2$, then $b_3=a_3$. Moreover, 
	\begin{align*}
		E(G') \supseteq &~\{xb_1,xa_2,xa_3,a_1a_2,a_1a_3,a_2a_3\} := F \\ 
		E(G') \subseteq &~F \cup \{b_1a_2,b_1a_3\}.
	\end{align*}		
	\item If $b_1 \neq a_1$, $b_2\neq a_2$ and $b_3=a_3$, then
	\begin{align*}
	E(G') \supseteq &~\{xb_1,xb_2,xa_3,a_1a_2,a_1a_3,a_2a_3\} := F \\ 
	E(G') \subseteq &~F \cup \{xa_2,b_1a_2,b_1b_2,b_1a_3,b_2a_1,b_2a_3\}.
	\end{align*}
	\item If $b_1 \neq a_1$, $b_2\neq a_2$ and $b_3\neq a_3$, then
	\begin{align*}
	E(G') \supseteq &~\{xb_1,xb_2,xb_3,a_1a_2,a_1a_3,a_2a_3\} := F \\ 
	E(G') \subseteq &~F \cup \{xa_2,xa_3,b_1a_2,b_1b_2,b_1a_3,b_1b_3,b_2a_1,b_2a_3,b_2b_3,b_3a_1,b_3a_2\}.
	\end{align*}
\end{enumerate}
\end{claim}
\begin{proof}
This follows readily from our assumption $t(b_3)<t(b_2)<t(b_1)$ with~Claim~\ref{thm:special-ordering} and~Claim~\ref{thm:x-has-only-one-neighbor}.
\end{proof}

In our algorithm, we do not only consider graphs, but rather tuples containing a graph together with its list of vertex labels and a linear vertex ordering.
The algorithm is split into three parts.
\begin{itemize}
	\item Algorithm~\ref{algo:init-algo2} initializes all relevant tuples according to Claim~\ref{thm:startgraphs}.
	\item Algorithm~\ref{algo:expand} is the main procedure, where a certain tuple is extended in all possible relevant ways. This corresponds to a labeling step in our tripod traversal algorithm.
	
		The list Act is the list of currently active vertices. In each step of the traversal algorithm, an active vertex is picked.
		Then one neighbor each from the two other tripod classes is added to the set of active vertices (unless they have been visited before). Consequently, Algorithm~\ref{algo:expand} adds up to two new vertices to the tripod that correspond to these two neighbors.
		According to Claim~\ref{thm:active-doesnt-matter}, it does not matter which active vertex is picked next.
	
	The list Ord maintained by the algorithm corresponds to the ordering proposed by Claim~\ref{thm:special-ordering}.
	Algorithm~\ref{algo:expand} implicitly enumerates all orderings that obey the properties of Claim~\ref{thm:special-ordering}.
	Whenever the ordering of the vertices of the partial tripod generated so far does not obey the properties listed in Claim~\ref{thm:special-ordering}, we may prune.
	\item Algorithm~\ref{algo:feasible} is a subroutine we use to prune tuples we do not need to consider. We call a tuple \emph{prunable} if Algorithm~\ref{algo:feasible} applied to it returns the value \emph{false}.
	
	The criteria for a prunable tuple we apply (in that order) are as follows:
	\begin{itemize}
		\item an induced $P_6$,
		\item the graph is not 3-colorable,
		\item one of the properties of Claim~\ref{thm:x-has-only-one-neighbor} is violated,
		\item the ordering Ord does not obey Claim~\ref{thm:special-ordering}, and
		\item the partial tripod enumerated so far cannot be extended to a 4-vertex-critical 1-vertex extension of a tripod.
	\end{itemize}
\end{itemize}
We now come to the correctness proof of these algorithms.

\begin{algorithm}[ht!]
\caption{Generate 4-critical $P_6$-free 1-vertex extension of a tripod}
\label{algo:init-algo2}
  \begin{algorithmic}[1]
	\STATE $V:=\{x,b_1,a_1,a_2,a_3\}$ ~~~~~~~~~~~~~~~~~~~~~~~~~~~~~~~~~~ // in this case, $b_2=a_2$ and $b_3=a_3$
	\STATE $E^{\mbox{\scriptsize must}}:=\{xb_1,xa_2,xa_3,a_1a_2,a_1a_3,a_2a_3\}$
	\STATE $E^{\mbox{\scriptsize may}}:=\{b_1a_2,b_1a_3\}$
	\STATE $\mbox{Ord}:=(a_3,a_2,a_1,b_1,x)$ and $\mbox{Act}:=\{b_1\}$
	\STATE $A_1:=\{a_1,b_1\}$, $A_2:=\{a_2\}$, and $A_3:=\{a_3\}$
	\FOR{each $E \subseteq E^{\mbox{\scriptsize must}} \cup E^{\mbox{\scriptsize may}}$ with $E^{\mbox{\scriptsize must}} \subseteq E$}
			\STATE\label{step:expand-a2=b2-a3=b3} Expand($G=(V,E),A_1,A_2,A_3,\mbox{Ord},\mbox{Act}$)
	\ENDFOR

	\STATE $V:=\{x,b_1,b_2,a_1,a_2,a_3\}$ ~~~~~~~~~~~~~~~~~~~~~~~~~~~~~~~~~ // in this case, $b_2\neq a_2$ and $b_3=a_3$
	\STATE $E^{\mbox{\scriptsize must}}:=\{xb_1,xb_2,xa_3,a_1a_2,a_1a_3,a_2a_3\}$
	\STATE $E^{\mbox{\scriptsize may}}:=\{xa_2,b_1b_2,b_1a_2,b_1a_3,b_2a_1,b_2a_3\}$
	\STATE $\mbox{Ord}:=(a_3,a_2,a_1,b_2,b_1,x)$ and $\mbox{Act}:=\{b_1,b_2\}$
	\STATE $A_1:=\{a_1,b_1\}$, $A_2:=\{a_2,b_2\}$, and $A_3:=\{a_3\}$
	\FOR{each $E \subseteq E^{\mbox{\scriptsize must}} \cup E^{\mbox{\scriptsize may}}$ with $E^{\mbox{\scriptsize must}} \subseteq E$}
			\STATE\label{step:expand-one-equal} Expand($G=(V,E),A_1,A_2,A_3,\mbox{Ord},\mbox{Act}$)
	\ENDFOR

	\STATE $V:=\{x,b_1,b_2,b_3,a_1,a_2,a_3\}$ ~~~~~~~~~~~~~~~~~~~~~~~~~~~~~~~~ // in this case, $b_2\neq a_2$ and $b_3 \neq a_3$
	\STATE $E^{\mbox{\scriptsize must}}:=\{xb_1,xb_2,xb_3,a_1a_2,a_1a_3,a_2a_3\}$
	\STATE $E^{\mbox{\scriptsize may}}:=\{xa_2,xa_3,b_1b_2,b_1b_3,b_2b_3,b_1a_2,b_1a_3,b_2a_1,b_2a_3,b_3a_1,b_3a_2\}$
	\STATE $\mbox{Ord}:=(a_3,a_2,a_1,b_3,b_2,b_1,x)$ and $\mbox{Act}:=\{b_1,b_2,b_3\}$
	\STATE $A_1:=\{a_1,b_1\}$, $A_2:=\{a_2,b_2\}$, and $A_3:=\{a_3,b_3\}$
	\FOR{each $E \subseteq E^{\mbox{\scriptsize must}} \cup E^{\mbox{\scriptsize may}}$ with $E^{\mbox{\scriptsize must}} \subseteq E$}
			\STATE\label{step:expand-all-different} Expand($G=(V,E),A_1,A_2,A_3,\mbox{Ord},\mbox{Act}$)
	\ENDFOR
  \end{algorithmic}
\end{algorithm}

\begin{algorithm}[ht!]
\caption{Expand(Graph $G=(V,E)$, Set $A_1$, Set $A_2$, Set $A_3$, List Ord, Set Act)}
\label{algo:expand}
  \begin{algorithmic}[1]
		\IF{not Feasible($G$, $A_1$, $A_2$, $A_3$, Ord, Act)} \label{line:prunable}
			\STATE \textbf{return}
		\ENDIF
		
		\STATE\label{step:choose-active-vertex} pick a vertex $u$ from the set Act and let $\{i,j,k\}=\{1,2,3\}$ be such that $u \in A_i$
		\STATE let $u_j$ be the $<_{\mbox{\scriptsize Ord}}$-minimal neighbor of $u$ in $A_j$, if existent, and let $u_k$ be defined accordingly\\
		~~~// we write $u <_{\mbox{\scriptsize Ord}} v$ whenever $u$ appears before $v$ in the list Ord
		\STATE\label{step:two-new-guys} let $v_j,v_k$ be two entirely new vertices
		\FOR{all ways of inserting $v_j$ and $v_k$ into the list Ord such that\vspace{-0.2cm}
			\begin{enumerate}[(a)]
				\item $a_1 <_{\mbox{\scriptsize Ord}} v_j,v_k <_{\mbox{\scriptsize Ord}} u$,\vspace{-0.2cm}
				\item $v_j <_{\mbox{\scriptsize Ord}} u_j$, if existent, and $v_k <_{\mbox{\scriptsize Ord}} u_k$, if existent\vspace{-0.2cm}
			\end{enumerate}}\label{step:insert-two-guys}
			\STATE\label{step:E*-two-new-neighbors} \vspace{-0.2cm}
			\begin{align*}
				E^*:= {}&{} \{w v_j : w \in A_i \cup A_k \mbox{ is active}\} \cup \{w v_k : w \in A_i \cup A_j \mbox{ is active}\} \cup \{x v_j, x v_k\} \cup \{v_j,v_k\}\\\vspace{-0.2cm}
							{}&{} \cup \{w v_j : w \in A_i \cup A_k \mbox{ is inactive and has a neighbor $w' \in A_j$ with $w' <_{\mbox{\scriptsize Ord}} v_j$}\}\\\vspace{-0.2cm}
							{}&{} \cup \{w v_k : w \in A_i \cup A_k \mbox{ is inactive and has a neighbor $w' \in A_k$ with $w' <_{\mbox{\scriptsize Ord}} v_k$}\}\vspace{-0.2cm}
			\end{align*}
			\FOR{all subsets $E'$ of $E^*$}
				\STATE\label{step:classes-Act-two-new-guys} $A_i':=A_i$, $A_j':=A_j\cup \{v_j\}$, $A_k':=A_k\cup \{v_k\}$, and $\mbox{Act}':=(\mbox{Act} \setminus \{u\}) \cup \{v_j,v_k\}$
				\STATE\label{step:Ord-two-new-guys} let $\mbox{Ord}'$ be Ord where $v_j$ and $v_k$ are inserted in the position we currently consider
				\STATE Expand($(V\cup \{v_j,v_k\},E\cup E'), A_1', A_2', A_3', \mbox{Ord}', \mbox{Act}')$
			\ENDFOR
		\ENDFOR
		
		\FOR{$r = j,k$}
			\IF{$u_r$ is existent and $u_r <_{\mbox{\scriptsize Ord}} u$}
				\STATE let $\{r,s\} = \{j,k\}$
				\FOR{all ways of inserting $v_s$ into the list Ord such that $a_1 <_{\mbox{\scriptsize Ord}} v_s <_{\mbox{\scriptsize Ord}} u$}
					\STATE \vspace{-0.2cm}
					\begin{align*}
					E^*:= {}&{} \{w v_s : w \in A_i \cup A_r \mbox{ is active}\} \cup \{x v_s\}\\
					{}&{} \cup \{w v_s : w \in A_i \cup A_r \mbox{ is inactive and has a neighbor $w' \in A_s$ with $w' <_{\mbox{\scriptsize Ord}} v_s$}\}\vspace{-0.2cm}
					\end{align*}
					\FOR{all subsets $E'$ of $E^*$}
						\STATE $A_i':=A_i$, $A_s':=A_s\cup \{v_s\}$, $A_r':=A_r$, and $\mbox{Act}':=(\mbox{Act} \setminus \{u\}) \cup \{v_s\}$
						\STATE let $\mbox{Ord}'$ be Ord where $v_s$ is inserted in the position we currently consider
						\STATE Expand($(V\cup \{v_s\},E\cup E'), A_1', A_2', A_3', \mbox{Ord}', \mbox{Act}')$
					\ENDFOR
				\ENDFOR
			\ENDIF
		\ENDFOR

		\IF{both $u_j$ and $u_k$ exist and $u_j,u_k <_{\mbox{\scriptsize Ord}} u$}
			\STATE Expand($G, A_1, A_2, A_3, \mbox{Ord}, \mbox{Act}\setminus \{u\})$
		\ENDIF
  \end{algorithmic}
\end{algorithm}

\begin{algorithm}[ht!]
\caption{Feasible(Graph $G=(V,E)$, Set $A_1$, Set $A_2$, Set $A_3$, List Ord, Set Act)}
\label{algo:feasible}
  \begin{algorithmic}[1]
		\IF{$G$ contains a $P_6$} \label{line:P6check}
			\STATE \textbf{return} \textbf{false}
		\ENDIF
		\IF{$G$ is not $3$-colorable}\label{line:4-chromatic-tripod}
			\STATE \textbf{return} \textbf{false}
		\ENDIF
		\IF{$x$ has at least two neighbors in $A_1$} \label{line:NeighborNumberCheck-A1}
			\STATE \textbf{return} \textbf{false}
		\ENDIF
		\IF{$x$ has at least two neighbors in $A_2$ and at least two neighbors in $A_3$} \label{line:NeighborNumberCheck}
			\STATE \textbf{return} \textbf{false}
		\ENDIF
		\FOR{any two distinct vertices $u,v \in (V\setminus \{x\})$ with $u<_{\mbox{\scriptsize Ord}}v$}\label{step:parent-pruning-begin}
			\IF{$u \notin \mbox{Act}$}
				\STATE let $\{i,j,k\}=\{1,2,3\}$ be such that $u \in A_i$
				\STATE let $u_j$ be the $<_{\mbox{\scriptsize Ord}}$-minimal neighbor of $u$ in $A_j$, and let $u_k$ be defined accordingly
				\STATE let $\{i',j',k'\}=\{1,2,3\}$ be such that $v \in A_{i'}$
				\STATE let $v_{j'}$ be the $<_{\mbox{\scriptsize Ord}}$-minimal neighbor of $v$ in $A_{j'}$, if existent, and let $v_{k'}$ be defined accordingly				
				\IF{the following hold:\vspace{-0.2cm}
				\begin{enumerate}[(a)]
					\item $\{u_j,u_k\} \not\subseteq \{a_1,a_2,a_3\}$,\vspace{-0.2cm}
					\item $v_{j'}$ and $v_{k'}$ both exist, and\vspace{-0.2cm}
					\item $v_{j'},v_{k'} <_{\mbox{\scriptsize Ord}} u_{r}$ for some $r \in \{j,k\}$\vspace{-0.2cm}
				\end{enumerate}} \label{cond:parent-pruning-if}
					\STATE \textbf{return} \textbf{false} 	
				\ENDIF
			\ENDIF
		\ENDFOR\label{step:parent-pruning-end}

		\FOR{each $u \in (V\setminus \{x\})$} \label{step:minimality-pruning-begin}
			\STATE $W := \{v \in V : v <_{\mbox{\scriptsize Ord}} u\}$
			\STATE $B_i:=A_i \cap W$ for each $i=1,2,3$
			\WHILE{there is a vertex $v \in V \setminus (B_1 \cup B_2 \cup B_3 \cup \{u\})$ with neighbors in at least two of $B_1,B_2,B_3$}
				\IF{$v$ has neighbors in all three of $B_1,B_2,B_3$}
					\STATE \textbf{return} \textbf{false}
				\ELSE
					\STATE $B_i:=B_i \cup \{v\}$, where $B_i$ is the set that $v$ does not have neighbors in
				\ENDIF
			\ENDWHILE
		\ENDFOR\label{step:minimality-pruning-end}
		\STATE \textbf{return} \textbf{true}
  \end{algorithmic}
\end{algorithm}

\begin{lemma}\label{thm:tripod-generator-works}
Assume that Algorithm~\ref{algo:init-algo2} terminates and does never generate a tuple whose graph has $k+1$ or $k+2$ vertices, for some $k \ge 4$.
Then any 4-critical $P_6$-free graph which is a 1-vertex extension of a tripod has at most $k$ vertices.
\end{lemma}

To see this, let $G$ be a 4-critical $P_6$-free graph other than $K_4$ that is a 1-vertex extension of a tripod, with the notation from above.
We need the following claim.

\begin{claim}\label{thm:traversal-vs-generation}
There is a sequence of tuples $\Gamma^i=(G^i=(V^i,E^i),A_1^i,A_2^i,A_3^i,\mbox{Ord}^i,\mbox{Act}^i)$, $i=0,\ldots,r$, and a way of traversing the tripod $T$ in $r$ steps, in the way described above, for which the following holds, after possibly renaming vertices. Let $V(i)$ be set of all labeled vertices after the $i$-th iteration of the traversal, together with $x$, and let $\mbox{Act}(i)$ be the set of vertices which are active after the $i$-th iteration of the traversal, for $i=0,\ldots,r$.
\begin{enumerate}[(a)]
	\item\label{cond:properties-initial-call} At some point during the algorithm, Expand$(\Gamma^0)$ is called.
	\item During the procedure Expand$(\Gamma^i)$, $\Gamma^{i+1}$ is generated and so Expand$(\Gamma^{i+1})$ is called, for all $i=0,\ldots,r-1$.
	\item\label{cond:properties-tuple} The following holds, for all $i=0,\ldots,r$.
		\begin{enumerate}[(i)]
			\item $G|V(i) = G^i$, and in particular $A_j \cap V(i) = A_j^i$, for all $j=1,2,3$,
			\item $\mbox{Act}^i = \mbox{Act}(i)$, and
			\item\label{cond:t-is-Ord} for any two $u,v \in V(i)$ with $t(u)<t(v)$, $u <_{\mbox{\scriptsize Ord}^i} v$.
		\end{enumerate}
\end{enumerate}
\end{claim}
\begin{proof}
Since $G$ is not $K_4$ we may assume that $b_1 \neq a_1$, by Claim~\ref{thm:startgraphs}.

If $b_2 = a_2$, then Claim~\ref{thm:startgraphs} implies $b_3 = a_3$, and $\Gamma^0$ is generated by Algorithm~\ref{algo:init-algo2}.
Here, $\Gamma^0=(G^0=(V^0,E^0),A_1^0,A_2^0,A_3^0,\mbox{Ord}^0,\mbox{Act}^0)$ with
\begin{itemize}
	\item $V^0=\{a_1,b_1,a_2,a_3,x\}$ and $E^0=E(G|V^0)$,
	\item $A_1^0=\{a_1,b_1\}$, $A_2^0=\{a_2\}$, and $A_3^0=\{a_3\}$, and
	\item $\mbox{Ord}^0 = (a_3,a_2,a_1,b_1,x)$, and $\mbox{Act}^0 = \{b_1\}$.
\end{itemize}
By Claim~\ref{thm:startgraphs}, $E^{\mbox{\scriptsize must}} \subseteq E^0 \subseteq E^{\mbox{\scriptsize must}} \cup E^{\mbox{\scriptsize may}}$.
The cases when $a_2\neq b_2$ but $a_3=b_3$ resp.~$a_3 \neq b_3$ are dealt with similarly.
This proves~\eqref{cond:properties-tuple} for $i=0$.

For the inductive step assume that for some $s \in \{0,\ldots,r-1\}$ the tuple $\Gamma^s$ has the properties mentioned in~\eqref{cond:properties-tuple}.
We first prove that $\Gamma^{s+1}$ is generated while Expand$(\Gamma^s)$ is processed, and that $\Gamma^{s+1}$ has the properties mentioned in~\eqref{cond:properties-tuple}.

First we discuss why Algorithm~\ref{algo:feasible} returns \emph{true} on the input $\Gamma^s$.
Clearly $G^s=G|V(s) \neq G$ is 3-colorable and $P_6$-free, and so the if-conditions in lines~\ref{line:4-chromatic-tripod} and~\ref{line:P6check} both do not apply.
Also, the if-conditions in the lines~\ref{line:NeighborNumberCheck-A1} and~\ref{line:NeighborNumberCheck} does not apply to $\Gamma^s$ due to Claim~\ref{thm:x-has-only-one-neighbor} applied to $G$ together with~\eqref{cond:properties-tuple}.(i) in the case $i=s$.

During the steps~\ref{step:parent-pruning-begin}-\ref{step:parent-pruning-end}, the if-condition in line~\ref{cond:parent-pruning-if} never applies due to Claim~\ref{thm:special-ordering}.
To see this, pick two distinct vertices $u,v \in (V^s\setminus \{x\})$ with $u<v$ and $u \notin \mbox{Act}^s$.
Let $\{i,j,k\}=\{1,2,3\}$ be such that $u \in A_i^s$, let $u_j$ be the $<_{\mbox{\scriptsize Ord}^s}$-minimal neighbor of $u$ in $A_j^s$, and let $u_k$ be defined accordingly, let $\{i',j',k'\}=\{1,2,3\}$ be such that $v \in A_{i'}^s$, and let $v_{j'}$ be the $<_{\mbox{\scriptsize Ord}^s}$-minimal neighbor of $v$ in $A_{j'}^s$, if existent, and let $v_{k'}$ be defined accordingly.

Due to property~\eqref{cond:properties-tuple}.(iii), $t(u)<t(v)$.
Since $u \in V^s \setminus \mbox{Act}^s$, we know that $u \in V(s) \setminus \mbox{Act}(s)$, by~\eqref{cond:properties-tuple}.(i).
Thus, $n_j(u),n_k(u) \in V(s)$.
Moreover, by~\eqref{cond:properties-tuple}.(i), $n_j(u) = u_j$ and $n_k(u) = u_k$.
Now, if $v_{j'},v_{k'}$ both exist and $v_{j'},v_{k'} <_{\mbox{\scriptsize Ord}^s} u_{r}$ for some $r \in \{j,k\}$, then in particular $t(n_{j'}(v)),t(n_{k'}(v)) < t(n_r(u))$, in contradiction to~Claim~\ref{thm:special-ordering}.

Finally, $\Gamma^s$ is not pruned in the lines~\ref{step:minimality-pruning-begin}-\ref{step:minimality-pruning-end} since $G-u$ is 3-colorable for every $u \in V$.

Now we argue why $\Gamma^{s+1}$ is constructed and carries the desired properties.
If $s=0$, the case is clear, so we may assume that $s > 0$.
Say that, in the procedure Expand$(\Gamma^s)$, vertex $u$ is picked in line~\ref{step:choose-active-vertex} of Algorithm~\ref{algo:expand}.
Let us say that $u \in A_i^{s}$, where $\{i,j,k\}=\{1,2,3\}$.
In the traversal procedure, $n_j(u)$ and $n_k(u)$ are now visited and made active, in case they are not in $V(s)$ already.

Let us first assume that $n_j(u),n_k(u) \notin V(s)$, and let $v_j,v_k$ be the two entirely new vertices picked in line~\ref{step:two-new-guys}.
Due to the definition of tripods, $t(a_1) < t(n_j(u)),t(n_k(u)) < t(u)$, and $$t(n_\ell(u)) < \min(\{t(w):w \in N_G(u) \cap A_\ell\} \cup \{\infty\}) \mbox{ for } \ell=j,k.$$
Consequently, the algorithm considers in line~\ref{step:insert-two-guys} inserting the two new vertices $v_j$ and $v_k$ into Ord$^s$ such that~\eqref{cond:properties-tuple}.(iii) holds, where we identify $v_j$ with $n_j(u)$ and $v_k$ with $n_k(u)$.
Moreover, $E^*$ in line~\ref{step:E*-two-new-neighbors} contains all edges incident to $n_j(u)$ and $n_k(u)$ in $G|V(s)$, due to the definition of $n_j(u)$ and $n_k(u)$.
Due to steps~\ref{step:classes-Act-two-new-guys} and~\ref{step:Ord-two-new-guys}, the tuple $\Gamma^{s+1}$ is indeed generated, and Expand$(\Gamma^{s+1})$ is called, where
\begin{itemize}
	\item $G^{s+1}=G|(V(s)\cup\{n_j(u),n_k(u)\})=G|V(s+1)$, and in particular $A_i^{s+1}=A_i^{s}=V(s)\cap A_i = V(s+1)\cap A_i$, and $A_\ell^{s+1}=A_\ell^{s}\cup\{v_\ell=n_\ell(u)\}=V(s+1)\cap A_\ell$ for $\ell=j,k$,
	\item $\mbox{Act}^{s+1} = (\mbox{Act}^{s+1} \setminus \{u\}) \cup\{v_j,v_k\} =  (\mbox{Act}(s) \setminus \{u\}) \cup\{n_j(u),n_k(u)\} = \mbox{Act}(s+1)$, and
	\item for any two vertices $u,v \in V(s+1)$ with $t(u)<t(v)$, $u <_{\mbox{\scriptsize Ord}^{s+1}} v$.
\end{itemize}
The cases when $n_j(u)$ and/or $n_k(u)$ have been active before are handled analogously.
This completes the proof of~Claim~\ref{thm:traversal-vs-generation}.
\end{proof}

Next we derive Lemma~\ref{thm:tripod-generator-works}.

\begin{proof}[Proof of~Lemma~\ref{thm:tripod-generator-works}]
Like above, $\Gamma^r$ is not pruned in step~\ref{line:prunable} during the procedure of Expand$(\Gamma^r)$.
Since $G^r=G|V(r)=G$, $G$ is indeed generated by the algorithm.
As $|V(G^s)|+2 \ge |V(G^{s+1})|$ for all $s=0,\ldots,r-1$, $G$ has at most $k$ vertices.
\end{proof}

We implemented this set of algorithms in C with some further optimizations. 
A crucial detail is how the active vertex is picked in line~\ref{step:choose-active-vertex} of Algorithm~\ref{algo:expand}.
The following choice seemed to terminate most quickly.
\begin{itemize}
	\item If the graph which is currently expanded has at most 12 vertices, we pick the $\mbox{Ord}$-maximal active vertex in line~\ref{step:choose-active-vertex}.
	\item If the graph has more than 12 vertices, we pick the active vertex for which the number of non-prunable tuples generated from it is minimum. This is done by trying to extend every active vertex once without iterating any further and counting the number of non-prunable tuples generated.
\end{itemize}
With this choice, our program does indeed terminate (in about 60 hours) and the largest non-prunable generated graph has $18$ vertices.
Together with Lemma~\ref{thm:tripod-generator-works}, we arrive at Lemma~\ref{thm:uncontractingK4}. 
Table~\ref{table:counts-generated-tripods} shows the number of non-prunable tuples generated by the program.


\begin{table}[ht!]
\centering
\begin{tabular}{| l || c c c c c |}
\hline 
$|V(G)|$            		&       5 &         6 &         7 & 		     8 &      9 \\
\# non-prunable tuples  &       3 &        67 &     2,010 &     11,726 & 81,523 \\
\hline
$|V(G)|$            		&      10 &			   11 & 			 12 &  				13 & 				 14 \\
\# non-prunable tuples  & 388,190 & 1,234,842 & 3,380,785 & 10,669,960 &  16,322,798 \\
\hline
$|V(G)|$            		&      15 &        16 &        17 &         18 &  19, 20 \\
\# non-prunable tuples  & 137,031 &    49,506 &     2,865 &        330 &       0 \\
\hline
\end{tabular}
\caption{Counts of the number of non-prunable tuples generated by our implementation of Algorithm~\ref{algo:init-algo2}}
\label{table:counts-generated-tripods}
\end{table}

In order to be sure the algorithm is implemented correctly, we also modified the program so it collects all 4-critical graphs found along the way, similar to line~\ref{line:k-critical} of Algorithm~\ref{algo:construct}.
As expected, all 4-critical $P_6$-free 1-vertex extensions of a tripod in $\mathcal L$ were found. In the Appendix we describe in more detail how we tested the correctness of our implementation and the source code of the program can be downloaded from~\cite{tripodgenerator-site}.

\section{Obstructions up to 28 vertices}\label{sec:smallgraphs}

In this section we prove the following result.

\begin{lemma}\label{thm:small-obstructions}
Let $G$ be a 4-critical $P_6$-free graph.
If $|V(G)| \le 28$, then $G$ is contained in $\mathcal L$.
\end{lemma}

For the proof of this result, we run the enumeration algorithm of Section~\ref{algo:init-algo}, with the following modifications.
In line~\ref{line:isocheck} of Algorithm~\ref{algo:construct}, we do not discard a graph if it contains a diamond, only when it is not $P_6$-free.
Moreover, we discard a graph if it contains more than 28 vertices.
This procedure terminates exactly with the list $\mathcal L$ (note that the largest graph in $\mathcal L$ has 16 vertices). 
Table~\ref{table:counts_P6free} shows the number of graphs generated by the algorithm on each relevant number of vertices. This computation took approximately 9 CPU years on a cluster.

\begin{table}[ht!]
\centering
\begin{tabular}{| l || c c c c c c |}
\hline 
$|V(G)|$            & 	 	    5 &         6 &         7 &         8 &     	     9 &        10 \\
\# graphs generated &         1 &         7 &        45 &       253 &   	   1,385 &     5,402 \\
\hline
$|V(G)|$            &        11 &        12 &        13 &        14 & 	        15 &        16 \\
\# graphs generated &    12,829 &     24,802 &   36,435 &    41,422 & 	    42,769 &    46,176 \\
\hline
$|V(G)|$            &        17 &        18 &        19 &        20 &  			    21 &      22 \\
\# graphs generated &    54,001 &    70,205 &    99,680 &   145,968 &			 233,687 & 382,762 \\
\hline
$|V(G)|$            &        23 &        24 &        25 &        26 &   	      27 & 28\\
\# graphs generated &   696,462 & 1,430,280 & 3,002,407 & 6,410,184 & 	13,703,206 &  30,764,536\\
\hline
\end{tabular}
\caption{Counts of the number of $P_6$-free graphs generated by our implementation of Algorithm~\ref{algo:init-algo} without testing for induced diamonds}
\label{table:counts_P6free}
\end{table}

\section{Proof of Theorem~\ref{thm:N3P6}}\label{sec:proof-of-main-result}

Let $G$ be a 4-critical $P_6$-free graph.
Suppose that $G \notin \mathcal L$ and that $|V(G)|$ is minimal with respect to this property.
If $G$ is diamond-free, Lemma~\ref{lem:diamondfree} implies $G\in \mathcal L$, a contradiction.
We may thus assume that there is a maximal tripod $T=(A_1,A_2,A_3)$ in $G$ which is not just a triangle.

Suppose that there is some vertex $x \in V(G) \setminus V(T)$ with a neighbor in each $A_i$, $i=1,2,3$.
Then $V(G)=V(T) \cup \{x\}$, and so $|V(G)| \le 18$ by Lemma~\ref{thm:uncontractingK4}.
By Lemma~\ref{thm:small-obstructions}, $G \in \mathcal L$, a contradiction.

So, we may assume that no vertex has a neighbor in all three classes of $T$.
Let $G'$ be the graph obtained by contracting $T$ in $G$.
By Lemma~\ref{lem:contraction-is-safe} we know that $G'$ is $P_6$-free and not 3-colorable.
We may thus pick a 4-critical $P_6$-free subgraph $H$ of $G'$.

Since $G$ was chosen to have a minimal number of vertices among all $4$-critical $P_6$-free graphs not in $\mathcal L$, we may assume that $H \in \mathcal L$.
Let $T'$ be the triangle in $G'$ obtained by contracting $T$.
If $H$ is a $K_4$, $|V(H)\cap V(T')| \le 2$, since no vertex of $G$ has a neighbor in all three classes of $T$. Since $H \in \mathcal{L}$, 
$|V(H)| \leq  16$.
Now Lemma~\ref{thm:10} and the minimality of $G$  imply that  
$|V(G)|\leq |V(H)|+12 \leq 28$. Consequently, by  Lemma~\ref{thm:small-obstructions}, 
$G \in \mathcal L$.

\section{$P_7$-free obstructions}\label{sec:P7}

This section is devoted to the following unpublished observation by Pokrovskiy~\cite{Pok14}. 

\begin{lemma}\label{thm:construction}
There are infinitely many 4-critical $P_7$-free graphs.
\end{lemma}

In the proof we construct an infinite family of 4-vertex-critical $P_7$-free graphs, i.e., $P_7$-free graphs which are 4-chromatic but every proper induced subgraph is 3-colorable.
This means that there is also an infinite number of $4$-critical $P_7$-free graphs.
Note that, indeed, not all members of our family are 4-critical $P_7$-free.

\begin{figure}
	\centering
		\includegraphics[width=0.40\textwidth]{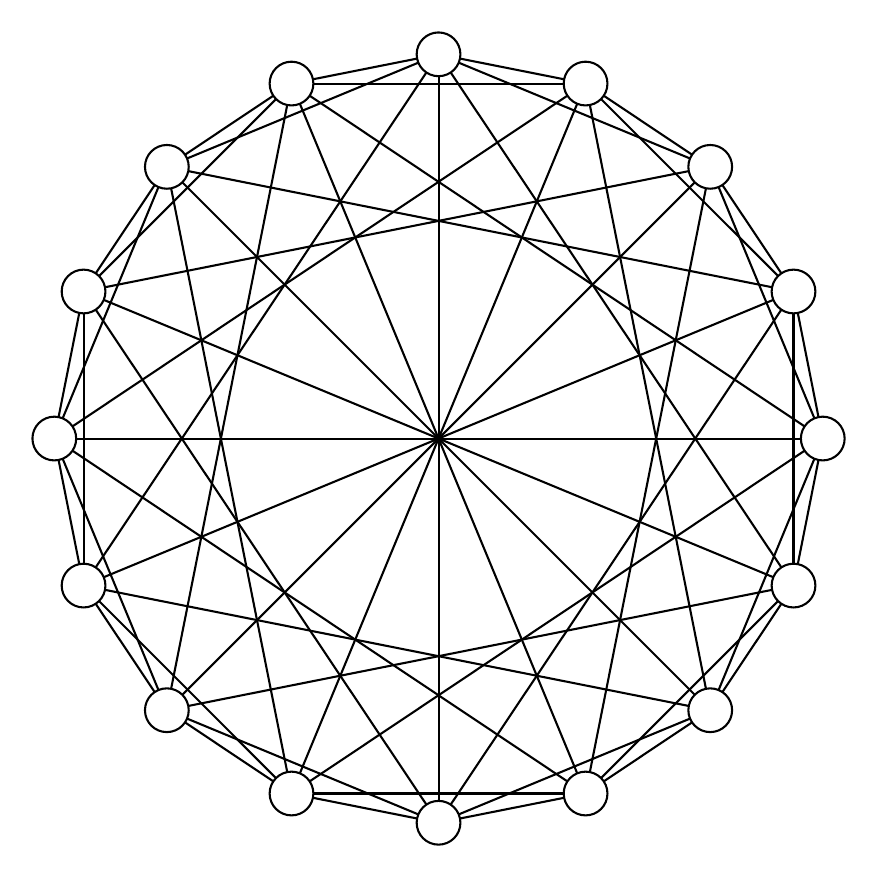}
	\label{fig:construction}
	\caption{A circular drawing of $G_{5}$}
	\label{fig:G16}
\end{figure}

\begin{proof}[Proof of Lemma~\ref{thm:construction}]
Consider the following construction.
For each $r \ge 1$, $G_r$ is a graph defined on the vertex set $v_0,\ldots,v_{3r}$.
The graph $G_5$ is shown in Fig.~\ref{fig:G16}. 
A vertex $v_i$, where $i \in \{0,1,\ldots,3r\}$, is adjacent to $v_{i-1}$, $v_{i+1}$, and $v_{i+3j+2}$, for all $j \in \{0,1,\ldots,r-1\}$.
Here and throughout the proof, we consider the indices to be taken modulo $3r+1$.

First we observe that, up to permuting the colors, there is exactly one 3-coloring of $G_r-v_0$.
Indeed, we may w.l.o.g.~assume that $v_i$ receives color $i$, for $i=1,2,3$, since $\{v_1,v_2,v_3\}$ forms a triangle in $G_r$.
Similarly, $v_4$ receives color 1, $v_5$ receives color 2 and so on.
Finally, $v_{3r}$ receives color 3.
Since the coloring was forced, our claim is proven.

In particular, $G_r$ is not 3-colorable, since $v_0$ is adjacent to all of $v_1,v_2,v_{3r}$.
As the choice of $v_0$ was arbitrary, we know that $G_r$ is 4-vertex-critical.

It remains to prove that $G_r$ is $P_7$-free.
Suppose that $P=x_1$-$x_2$-$\ldots$-$x_7$ is an induced $P_7$ in $G_r$.
To simplify the argumentation, we assume $G_r$ to be equipped with the proper coloring described above.
That is, $v_0$ has color 4, and, for all $i=0,\ldots,r-1$ and $j=1,2,3$, the vertex $v_{3i+j}$ is colored with color $j$.
Let $X_i$ denote the set of vertices of color $i$, for $i=1,2,3,4$.

If $r \le 2$, $|V(G_r)| \le 7$, and so we are done since obviously $G_2$ is not isomorphic to $P_7$.
Therefore, we may assume $r \ge 3$ and, since $G_r$ is vertex-transitive, w.l.o.g.~$v_0 \notin V(P)$.
Hence, $P$ is an induced $P_7$ in the graph $H:=G_r-v_0$, which we consider from now on.

First we suppose that some vertices of $P$ appear consecutively in the ordering $v_1,\ldots,v_{3r}$.
That is, w.l.o.g.~$x_i = v_j$ and $x_{i+1} = v_{j+1}$ for some $i \in \{1,\ldots,6\}$ and $j \in \{1,\ldots,3r-1\}$.
Since $P$ is an induced path, we know that neither of $v_{j-1}$ and $v_{j+2}$, if existent, are contained in $P$.
Thus, we may assume that $j=1$, and so $v_3 \notin V(P)$.
Recall that $N_H(v_1)\setminus \{v_2\} = X_3$ and $N_H(v_2)\setminus \{v_3\} = X_1$.
Thus, $|N_H(x_i) \cap V(P)|\le 2$ implies $|X_3 \cap V(P)|\le 1$, and similarly $|N_H(x_{i+1}) \cap V(P)|\le 2$ implies $|X_1 \cap V(P)|\le 2$.
Therefore, $|X_2 \cap V(P)| = 4$, which means that $x_1,x_3,x_5,x_7 \in X_2$.
But this is a contradiction to the fact that $N_H(v_1)\setminus \{v_2\} = X_3$.

Hence, no two vertices of $P$ appear consecutively in the ordering $v_1,\ldots,v_{3r}$.
For simplicity, let us say that a vertex $v_i$ is \emph{left of} (\emph{right of}) a vertex $v_j$ if $i < j$ (if $i > j$).
We now know the following.
Let $x \in V(P)$ be left of $y \in V(P)$. Then $xy \in E$ if and only if $x \in X_1$ and $y \in X_3$, $x \in X_2$ and $y \in X_1$, or $x \in X_3$ and $y \in X_2$.
Below we make frequent use of this fact without further reference.

W.l.o.g.~$x_1 \in X_1$ and $x_2 \in X_2$.
In particular, $x_2$ is left of $x_1$.
We now distinguish the possible colorings of the remaining vertices of $P$, obtaining a contradiction in each case.

\textbf{Case 1.} \emph{$x_3 \in X_1$.}
	
In this case, $x_3$ must be right of $x_2$.

\textbf{Case 1.1.} \emph{$x_4 \in X_2$.}

In this case, $x_4$ is right of $x_1$, and in turn $x_3$ is right of $x_4$.
Hence, $x_5$ cannot be in $X_1$, since then it must be right of $x_4$ but left of $x_2$.
So, $x_5 \in X_3$, and thus $x_5$ is between $x_2$ and $x_1$.

\textbf{Case 1.1.1.} \emph{$x_6 \in X_1$.}

Then $x_6$ must be left of $x_2$.
If $x_7 \in X_2$, it must be left of $x_6$ but right of $x_3$, a contradiction.
Otherwise if $x_7 \in X_3$, it must be left of $x_1$ but right of $x_4$, another contradiction.

\textbf{Case 1.1.2.} \emph{$x_6 \in X_2$.}

In this case $x_6$ must be right of $x_3$.
If $x_7 \in X_1$, it must be left of $x_4$ but right of $x_6$, a contradiction.
Otherwise if $x_7 \in X_3$, it must be left of $x_1$ but right of $x_4$, again a contradiction.

\textbf{Case 1.2.} \emph{$x_4 \in X_3$.}

In this case, $x_4$ is right of $x_3$, and in turn $x_1$ is right of $x_4$.

\textbf{Case 1.2.1.} \emph{$x_5 \in X_1$.}

So, $x_5$ must be left of $x_2$.
Hence, $x_6$ cannot be in $X_2$, since then $x_6$ must be left of $x_5$ and right of $x_1$.
Thus, $x_6 \in X_3$, which means that $x_6$ is between $x_2$ and $x_3$.

If $x_7 \in X_1$, it must be left of $x_6$ but right of $x_1$, a contradiction.
Otherwise if $x_7 \in X_2$, it must be left of $x_4$ but right of $x_1$, another contradiction.

\textbf{Case 1.2.2.} \emph{$x_5 \in X_2$.}

So, $x_5$ must be right of $x_1$.
Clearly $x_6 \notin X_1$, for then it must be left of $x_2$ but right of $x_5$.
So, $x_6 \in X_3$, and thus $x_6$ is between $x_2$ and $x_3$.

If $x_7 \in X_1$, it must be left of $x_6$ but right of $x_4$, a contradiction.
Otherwise if $x_7 \in X_2$, it must be left of $x_4$ but right of $x_1$, another contradiction.

\textbf{Case 2.} \emph{$x_3 \in X_3$.}
	
In this case, $x_3$ must be left of $x_2$.

\textbf{Case 2.1.} \emph{$x_4 \in X_1$.}

Then $x_4$ is left of $x_3$, and thus also $x_1$ and $x_2$.

If $x_5 \in X_2$, $x_5$ must be left of $x_4$ but right of $x_1$, a contradiction.
So, $x_5 \in X_3$.
Then $x_5$ must be between $x_2$ and $x_1$.
If $x_6 \in X_2$, it must be right of $x_5$ but left of $x_3$, a contradiction.
So, $x_6 \in X_1$, and thus $x_6$ must be between $x_3$ and $x_2$.

If $x_7 \in X_2$, it must be left of $x_6$ but right of $x_1$, a contradiction.
Hence, $x_7 \in X_3$.
But now $x_7$ must be right of $x_6$ and left of $x_4$, another contradiction.

\textbf{Case 2.2.} \emph{$x_4 \in X_2$.}

Then $x_4$ must be right of $x_1$.

If $x_5 \in X_1$, it must be right of $x_4$ but left of $x_2$, a contradiction.
So, $x_5 \in X_3$, and thus $x_5$ is between $x_2$ and $x_1$.

If $x_6 \in X_1$, it must be between $x_3$ and $x_2$.
If, moreover, $x_7 \in X_2$, $x_7$ is left of $x_6$ but right of $x_1$, a contradiction.
Similarly, if $x_7 \in X_3$, $x_7$ is left of $x_1$ but right of $x_4$, another contradiction.

We thus know $x_6 \in X_2$.
But then $x_6$ must be left of $x_3$ and right of $x_5$, a contradiction.

Summing up, $G_r$ is $P_7$-free, and this completes the proof.
\end{proof}

We also modified Algorithm~\ref{algo:construct} to generate 4-critical $P_7$-free graphs. 
As one would expect, the number of obstructions is much larger than in the $P_6$-free case.
Table~\ref{table:counts_N3P7} contains the counts of all 4-critical and 4-vertex-critical $P_7$-free graphs up to 15 vertices. 

\begin{table}[ht!]
\centering
\begin{tabular}{| c | r | r |}
\hline 
Vertices & Critical graphs & Vertex-critical graphs\\
\hline 
4  &  1  &  1\\
6  &  1  &  1\\
7  &  2  &  7\\
8  &  5  &  8\\
9  &  21  &  124\\
10  &  99  &  2,263\\
11  &  212  &  1,771\\
12  &  522  &  6,293\\
13  &  679  &  15,064\\
14  &  368  &  4,521\\
15  &  304  &  2,914\\
\hline
$\le 15$  &  2,214  &  32,967\\
\hline
\end{tabular}
\caption{Counts of all 4-critical and 4-vertex-critical $P_7$-free graphs up to 15 vertices}

\label{table:counts_N3P7}
\end{table}


\subsection*{Acknowledgments}

We thank Alexey Pokrovskiy for pointing us to the existence of an infinite family of $P_7$-free obstructions.

Several of the computations for this work were carried out using the Stevin Supercomputer Infrastructure at Ghent University.


\bibliographystyle{amsplain}
\bibliography{bnm-j,references,bnm}

\providecommand{\bysame}{\leavevmode\hbox to3em{\hrulefill}\thinspace}
\providecommand{\MR}{\relax\ifhmode\unskip\space\fi MR }
\providecommand{\MRhref}[2]{%
  \href{http://www.ams.org/mathscinet-getitem?mr=#1}{#2}
}
\providecommand{\href}[2]{#2}
\begin{thebibliography}{10}

\bibitem{BCMSSZ15}
F.~Bonomo, M.~Chudnovsky, P.~Maceli, O.~Schaudt, M.~Stein, and M.~Zhong,
  \emph{Three-coloring and list three-coloring graphs without induced paths on
  seven vertices}, to appear in Combinatorica,
  \url{http://doi.org/10.1007/s00493-017-3553-8}, 2017.

\bibitem{hog}
G.~Brinkmann, K.~Coolsaet, J.~Goedgebeur, and H.~M{\'e}lot, \emph{{House of
  Graphs: a database of interesting graphs}}, Discrete Applied Mathematics
  \textbf{161} (2013), no.~1-2, 311--314, Available at
  \url{http://hog.grinvin.org/}.

\bibitem{BHS09}
D.~Bruce, C.T. Ho{\`a}ng, and J.~Sawada, \emph{A certifying algorithm for
  3-colorability of ${P}_5$-free graphs}, Proceedings of the 20th International
  Symposium on Algorithms and Computation, Springer-Verlag, 2009, pp.~594--604.

\bibitem{chudnovsky_06}
M.~Chudnovsky, N.~Robertson, P.~Seymour, and R.~Thomas, \emph{{The Strong
  Perfect Graph Theorem}}, {Annals of Mathematics} \textbf{{164}} ({2006}),
  no.~{1}, {51--229}.

\bibitem{DBLP:conf/soda/ChudnovskyGSZ16}
Maria Chudnovsky, Jan Goedgebeur, Oliver Schaudt, and Mingxian Zhong,
  \emph{Obstructions for three-coloring graphs with one forbidden induced
  subgraph}, Proceedings of the Twenty-Seventh Annual {ACM-SIAM} Symposium on
  Discrete Algorithms, {SODA} 2016, Arlington, VA, USA, January 10-12, 2016,
  2016, pp.~1774--1783.

\bibitem{Erd59}
P.~Erd\H{o}s, \emph{Graph theory and probability}, Canadian Journal of
  Mathematics \textbf{{11}} (1959), {34--38}.

\bibitem{criticalpfree-site}
J.~Goedgebeur, Homepage of generator for 4-critical $P_t$-free graphs:
  \url{http://caagt.ugent.be/criticalpfree/}.

\bibitem{tripodgenerator-site}
\bysame, Homepage of generator for 4-critical $P_t$-free 1-vertex extensions of
  tripods: \url{http://caagt.ugent.be/tripods/}.

\bibitem{paper-other-cases}
J.~Goedgebeur and O.~Schaudt, \emph{{Exhaustive generation of $k$-critical
  $\mathcal H$-free graphs}}, Journal of Graph Theory \textbf{87(2)} (2018),
  188--207.

\bibitem{GJPS15}
P.A. Golovach, M.~Johnson, D.~Paulusma, and J.~Song, \emph{A survey on the
  computational complexity of coloring graphs with forbidden subgraphs},
  Journal of Graph Theory \textbf{84(4)} (2017), 331--363.

\bibitem{hell_14}
P.~Hell and S.~Huang, \emph{Complexity of coloring graphs without paths and
  cycles}, LATIN 2014: Theoretical Informatics, Springer, 2014, pp.~538--549.

\bibitem{HMRSV15}
C.T. Ho\`{a}ng, B.~Moore, D.~Recoskie, J.~Sawada, and M.~Vatshelle,
  \emph{Constructions of $k$-critical ${P}_5$-free graphs}, Discrete Applied
  Mathematics \textbf{{182}} (2015), {91--98}.

\bibitem{LU95}
F.~Lazebnik and V.A. Ustimenko, \emph{{Explicit construction of graphs with an
  arbitrary large girth and of large size}}, Discrete Applied Mathematics
  \textbf{{60}} ({1995}), {275--284}.

\bibitem{nauty-website}
B.D. McKay, \emph{{nauty User's Guide (Version 2.5)}}, Technical Report
  TR-CS-90-02, Department of Computer Science, Australian National University.
  The latest version of the software is available at
  \url{http://cs.anu.edu.au/~bdm/nauty}.

\bibitem{mckay_14}
B.D. McKay and A.~Piperno, \emph{Practical graph isomorphism, {II}}, Journal of
  Symbolic Computation \textbf{60} (2014), 94--112.

\bibitem{Pok14}
A.~Pokrovskiy, \emph{private communication}.

\bibitem{RST02}
B.~{Randerath}, I.~{Schiermeyer}, and M.~{Tewes}, \emph{Three-colorability and
  forbidden subgraphs. {I}{I}: polynomial algorithms}, Discrete Mathematics
  \textbf{251} (2002), 137--153.

\bibitem{SeyPriv}
P.~Seymour, \emph{Barbados workshop on graph coloring and structure, 2014}.

\bibitem{OEIS}
N.~Sloane, The on-line encyclopedia of integer sequences:
  \url{http://oeis.org/}.

\end{thebibliography}


\section*{Appendix 1: Correctness testing}

Since several results obtained in this paper rely on computations,
it is very important that the correctness of our programs has
been thoroughly verified to minimize the chance of programming errors. 
In the following subsections we explain how we tested the correctness of our implementations.

Since all of our consistency tests passed, we believe that this is strong evidence for the correctness of our implementations.

\subsection*{Appendix 1.1: Correctness testing of critical $P_t$-free graph generator}
We performed the following consistency tests to verify the correctness of our generator for $k$-critical $P_t$-free graphs (i.e.\ Algorithm~\ref{algo:init-algo}). The source code of this program can be downloaded from~\cite{criticalpfree-site}.

\begin{itemize}

\item We applied the program to generate critical graphs for cases which were already settled before in the literature and verified that our program indeed obtained the same results. More specifically we verified that our program yielded exactly the same results in the following cases:

\begin{itemize}
\item There are six 4-critical $P_5$-free graphs~\cite{BHS09}.
\item There are eight 5-critical $(P_5,C_5)$-free graphs~\cite{HMRSV15}.
\item The Gr\"otzsch graph is the only 4-critical $(P_6,C_3)$-free graph~\cite{RST02}.
\item There are four 4-critical $(P_6,C_4)$-free graphs~\cite{hell_14}.
\end{itemize}

\item We developed an independent generator for  $k$-critical $P_t$-free graphs by starting from the program \verb|geng|~\cite{nauty-website,mckay_14} (which is a generator for all graphs) and adding pruning routines to it for colorability and $P_t$-freeness. This generator cannot terminate, but we were able to independently verify the following results with it:

\begin{itemize}
\item We executed this program to generate all 4-critical $(P_6,\mbox{diamond})$-free graphs up to 16 vertices and it indeed yielded the same 6 critical graphs from~Lemma~\ref{lem:diamondfree}.
\item We executed this program to generate all 4-critical and 4-vertex-critical $P_6$-free graphs up to 16 vertices and it indeed yielded the same graphs from~Theorem~\ref{thm:N3P6} and Table~\ref{table:counts_animals}.

\item We executed this program to generate all 4-critical and 4-vertex-critical $P_7$-free graphs up to 13 vertices and it indeed yielded the same graphs from Table~\ref{table:counts_N3P7}.

\end{itemize}

\item We modified our program to generate all $P_t$-free graphs and compared it with the known counts of $P_t$-free graphs for $t=4,5$ on the On-Line Encyclopedia of Integer Sequences~\cite{OEIS} (i.e.\ sequences A000669 and A078564).

\item We modified our program to generate all $k$-colorable graphs and compared it with the known counts of  $k$-colorable graphs for $k=3,4$ on the On-Line Encyclopedia of Integer Sequences~\cite{OEIS} (i.e.\ sequences A076322 and A076323).

\item We determined all $k$-vertex-critical graphs in two independent ways and both methods yielded exactly the same results:

\begin{enumerate}
\item By modifying line~\ref{line:k-critical} of Algorithm~\ref{algo:construct} so it tests for $k$-vertex-criticality instead of $k$-criticality.

\item By recursively adding edges in all possible ways to the set of critical graphs (as long as the graphs remain $k$-vertex-critical) and  testing if the resulting graphs are $P_t$-free. 

\end{enumerate}

\end{itemize}

\subsection*{Appendix 1.2: Correctness testing of tripod generator}

We performed the following consistency tests to verify the correctness of our generator for 4-critical $P_6$-free 1-vertex extensions of tripods (i.e.\ Algorithm~\ref{algo:init-algo2}). The source code of this program can be downloaded from~\cite{tripodgenerator-site}.

\begin{itemize}
\item We wrote a program to test if a graph is a 1-vertex extension of a tripod and applied it to the 24 4-critical $P_6$-free graphs from Theorem~\ref{thm:N3P6}. 11 of those graphs are 1-vertex extensions of a tripod (i.e.\ $F_1$, $F_2$, $F_4$,  $F_6$, $F_7$, $F_9$, $F_{10}$, $F_{17}$, $F_{21}$, $F_{22}$ and $F_{23}$). We verified that our implementation of Algorithm~\ref{algo:init-algo2} indeed yielded exactly those 11 graphs which are a 1-vertex extension of a tripod (except $K_4$).

\item We used Algorithm~\ref{algo:init-algo} to generate all 4-critical $P_7$-free graphs up to 14 vertices. There are 1910 such graphs and 595 of them are 1-vertex extensions of a tripod (see Table~\ref{table:counts_4-crit_P7-free} for details). We modified our implementation of Algorithm~\ref{algo:init-algo2} to generate 4-critical $P_7$-free 1-vertex extensions of tripods and executed it up to 14 vertices. We verified that this indeed yields exactly those 595 graphs which are a 1-vertex extension of a tripod (except $K_4$).
\end{itemize}

\begin{table}[ht!]
\centering
\begin{tabular}{| c | r | r |}
\hline 
Vertices & Critical graphs & 1-vertex extensions\\
\hline 
4  &  1  &  1\\
6  &  1  &  1\\
7  &  2  &  1\\
8  &  5  &  4\\
9  &  21  &  14\\
10  &  99  &  56\\
11  &  212  &  87\\
12  &  522  &  141\\
13  &  679  &  196\\
14  &  368  &  94\\
\hline
$\le 14$  &  1,910  &  595\\
\hline
\end{tabular}
\caption{Counts of 4-critical $P_7$-free graphs up to 14 vertices and the number of those graphs which are 1-vertex extensions of a tripod}

\label{table:counts_4-crit_P7-free}
\end{table}

\section*{Appendix 2: Adjacency lists}

This section contains the adjacency lists of the 24 $4$-critical $P_6$-free graphs from Theorem~\ref{thm:N3P6}. The graphs are listed in the same order as in Fig.~\ref{fig:graphs1-13} and~\ref{fig:graphs14-24}.

\begin{itemize}
\item Graph $F_1$: \{0 : 1 2 3; 1 : 0 2 3; 2 : 0 1 3; 2 : 0 1 3\}
\item Graph $F_2$: \{0 : 2 3 5; 1 : 3 4 5; 2 : 0 4 5; 3 : 0 1 5; 4 : 1 2 5; 5 : 0 1 2 3 4\}
\item Graph $F_3$: \{0 : 2 4 5; 1 : 3 5 6; 2 : 0 4 6; 3 : 1 5 6; 4 : 0 2 6; 5 : 0 1 3; 6 : 1 2 3 4\}
\item Graph $F_4$: \{0 : 3 4 5; 1 : 3 5 6; 2 : 4 5 6; 3 : 0 1 4 6; 4 : 0 2 3 6; 5 : 0 1 2; 6 : 1 2 3 4\}
\item Graph $F_5$: \{0 : 3 4 5; 1 : 4 6 7; 2 : 5 6 7; 3 : 0 6 7; 4 : 0 1 5; 5 : 0 2 4; 6 : 1 2 3 7; 7 : 1 2 3 6\}
\item Graph $F_6$: \{0 : 3 5 6; 1 : 4 5 7; 2 : 5 6 7; 3 : 0 6 7; 4 : 1 6 7; 5 : 0 1 2; 6 : 0 2 3 4 7; 7 : 1 2 3 4 6\}
\item Graph $F_7$: \{0 : 3 4 5 7; 1 : 4 5 6; 2 : 5 6 7; 3 : 0 6 7; 4 : 0 1 7; 5 : 0 1 2; 6 : 1 2 3 7; 7 : 0 2 3 4 6\}
\item Graph $F_8$: \{0 : 3 5 7; 1 : 4 7 8; 2 : 5 6 7; 3 : 0 6 8; 4 : 1 7 8; 5 : 0 2 8; 6 : 2 3 8; 7 : 0 1 2 4; 8 : 1 3 4 5 6\}
\item Graph $F_9$: \{0 : 4 5 8; 1 : 4 7 8; 2 : 5 6 8; 3 : 6 7 8; 4 : 0 1 6 8; 5 : 0 2 7; 6 : 2 3 4 8; 7 : 1 3 5; 8 : 0 1 2 3 4 6\}
\item Graph $F_{10}$: \{0 : 4 5 7; 1 : 4 7 8; 2 : 5 6 7; 3 : 6 7 8; 4 : 0 1 6 8; 5 : 0 2 8; 6 : 2 3 4 8; 7 : 0 1 2 3; 8 : 1 3 4 5 6\}
\item Graph $F_{11}$: \{0 : 3 4 5 8; 1 : 4 5 6; 2 : 5 6 7 8; 3 : 0 6 7; 4 : 0 1 7 8; 5 : 0 1 2; 6 : 1 2 3 8; 7 : 2 3 4; 8 : 0 2 4 6\}
\item Graph $F_{12}$: \{0 : 3 6 9; 1 : 4 6 7; 2 : 5 7 8; 3 : 0 6 9; 4 : 1 8 9; 5 : 2 7 8; 6 : 0 1 3 8; 7 : 1 2 5 9; 8 : 2 4 5 6; 9 : 0 3 4 7\}
\item Graph $F_{13}$: \{0 : 4 6 9; 1 : 5 6 8; 2 : 6 8 9; 3 : 7 8 9; 4 : 0 7 8; 5 : 1 7 9; 6 : 0 1 2 7; 7 : 3 4 5 6; 8 : 1 2 3 4 9; 9 : 0 2 3 5 8\}
\item Graph $F_{14}$: \{0 : 4 5 7 9; 1 : 5 6 7; 2 : 6 7 8; 3 : 7 8 9; 4 : 0 6 8; 5 : 0 1 8 9; 6 : 1 2 4 9; 7 : 0 1 2 3; 8 : 2 3 4 5; 9 : 0 3 5 6\}
\item Graph $F_{15}$: \{0 : 4 5 8 9; 1 : 4 7 8 9; 2 : 5 6 8; 3 : 6 7 8; 4 : 0 1 6 9; 5 : 0 2 7; 6 : 2 3 4 9; 7 : 1 3 5; 8 : 0 1 2 3; 9 : 0 1 4 6\}
\item Graph $F_{16}$: \{0 : 5 6 7; 1 : 5 6 9; 2 : 5 8 9; 3 : 6 7 8; 4 : 7 8 9; 5 : 0 1 2 7 8; 6 : 0 1 3 8 9; 7 : 0 3 4 5 9; 8 : 2 3 4 5 6; 9 : 1 2 4 6 7\}
\item Graph $F_{17}$: \{0 : 3 5 6 9; 1 : 4 6 8; 2 : 5 6 7 8 9; 3 : 0 7 8 9; 4 : 1 7 9; 5 : 0 2 7 8; 6 : 0 1 2; 7 : 2 3 4 5; 8 : 1 2 3 5 9; 9 : 0 2 3 4 8\}
\item Graph $F_{18}$: \{0 : 5 6 10; 1 : 5 9 10; 2 : 6 7 10; 3 : 7 8 10; 4 : 8 9 10; 5 : 0 1 7 8; 6 : 0 2 8 9; 7 : 2 3 5 9; 8 : 3 4 5 6; 9 : 1 4 6 7; 10 : 0 1 2 3 4\}
\item Graph $F_{19}$: \{0 : 4 6 7 10; 1 : 5 9 10; 2 : 6 8 9 10; 3 : 7 8 9 10; 4 : 0 8 9; 5 : 1 9 10; 6 : 0 2 7; 7 : 0 3 6; 8 : 2 3 4 10; 9 : 1 2 3 4 5; 10 : 0 1 2 3 5 8\}
\item Graph $F_{20}$: \{0 : 5 10 11; 1 : 6 7 10 11; 2 : 6 9 10 11; 3 : 7 8 10 11; 4 : 8 9 10 11; 5 : 0 10 11; 6 : 1 2 8 10; 7 : 1 3 9; 8 : 3 4 6 11; 9 : 2 4 7; 10 : 0 1 2 3 4 5 6; 11 : 0 1 2 3 4 5 8\}
\item Graph $F_{21}$: \{0 : 4 6 7 10 11; 1 : 5 6 7 8 11; 2 : 6 8 10 11 12; 3 : 7 8 9 10 11 12; 4 : 0 8 9 12; 5 : 1 9 10 11 12; 6 : 0 1 2 9 10 12; 7 : 0 1 3 12; 8 : 1 2 3 4; 9 : 3 4 5 6 11; 10 : 0 2 3 5 6; 11 : 0 1 2 3 5 9; 12 : 2 3 4 5 6 7\}
\item Graph $F_{22}$: \{0 : 4 6 8 9 11 12; 1 : 5 6 7 10 11 12; 2 : 6 7 8 9 11 12; 3 : 9 10 11 12; 4 : 0 7 10 11 12; 5 : 1 8 9 12; 6 : 0 1 2 10; 7 : 1 2 4 9; 8 : 0 2 5 10 11; 9 : 0 2 3 5 7 10; 10 : 1 3 4 6 8 9; 11 : 0 1 2 3 4 8; 12 : 0 1 2 3 4 5\}
\item Graph $F_{23}$: \{0 : 4 6 7 9 10; 1 : 5 7 8 9; 2 : 6 7 9 10 11; 3 : 7 8 9 10 11 12; 4 : 0 8 9 10 11 12; 5 : 1 10 11 12; 6 : 0 2 8 11 12; 7 : 0 1 2 3 11 12; 8 : 1 3 4 6; 9 : 0 1 2 3 4 12; 10 : 0 2 3 4 5; 11 : 2 3 4 5 6 7; 12 : 3 4 5 6 7 9\}
\item Graph $F_{24}$: \{0 : 4 8 13 14 15; 1 : 5 8 10 14 15; 2 : 6 8 9 10 15; 3 : 7 8 9 10 11; 4 : 0 9 10 11 12; 5 : 1 9 11 12 13; 6 : 2 11 12 13 14; 7 : 3 12 13 14 15; 8 : 0 1 2 3 11 12 13; 9 : 2 3 4 5 13 14 15; 10 : 1 2 3 4 12 13 14; 11 : 3 4 5 6 8 14 15; 12 : 4 5 6 7 8 10 15; 13 : 0 5 6 7 8 9 10; 14 : 0 1 6 7 9 10 11; 15 : 0 1 2 7 9 11 12\}
\end{itemize}

\end{document}